\documentclass[11pt]{article}

\usepackage{amsmath, amsfonts, amssymb,latexsym,wasysym,graphicx,stmaryrd,tikz,txfonts}
\input epsf.sty

\newtheorem{Theorem}{Theorem}[section]
\newtheorem{Lemma}{Lemma}[section]
\newtheorem{Proposition}[Lemma]{Proposition}
\newtheorem{Corollary}[Lemma]{Corollary}
\newtheorem{Definition}[Lemma]{Definition}

\newcommand{\BEQ}{\begin{equation}}     
\newcommand{\BEA}{\begin{eqnarray}}
\newcommand{\BD}{\begin{displaymath}}
\newcommand{\EEQ}{\end{equation}}       
\newcommand{\EEA}{\end{eqnarray}}
\newcommand{\ED}{\end{displaymath}}

\newcommand{\del}{\delta}
\newcommand{\Del}{\Delta}
\newcommand{\eps}{\varepsilon}          

\newcommand{\Tr}{{\mathrm{Tr}}}

\newcommand{\R}{\mathbb{R}}

\newcommand{\Z}{\mathbb{Z}}

\newcommand{\eop}{\hfill $\Box$}        

\newcommand{\II}{{\rm i}}               
\newcommand{\half}{{1\over 2}}          

\catcode`\@=11
\def\numberbysection{\@addtoreset{equation}{section}
        \def\theequation{\thesection.\arabic{equation}}}
\numberbysection

\begin{document}

\vspace*{1.5cm}
\begin{center}
{\Large \bf Global existence for strong solutions of viscous Burgers equation. (1) The
bounded case}

\end{center}

\vspace{2mm}
\begin{center}
{\bf  J\'er\'emie Unterberger$^a$}
\end{center}

\vskip 0.5 cm
\centerline {$^a$Institut Elie Cartan,\footnote{Laboratoire 
associ\'e au CNRS UMR 7502} Universit\'e de Lorraine,} 
\centerline{ B.P. 239, 
F -- 54506 Vand{\oe}uvre-l\`es-Nancy Cedex, France}
\centerline{jeremie.unterberger@univ-lorraine.fr}

\vspace{2mm}
\begin{quote}

\renewcommand{\baselinestretch}{1.0}
\footnotesize
{We prove that the viscous Burgers equation $(\partial_t-\Del) u(t,x)+(u\cdot\nabla)u(t,x)=g(t,x),\ (t,x)\in\R_+\times\R^d$ $(d\ge 1)$
  has a globally defined
smooth solution in all dimensions provided the initial condition and the forcing term $g$ are smooth and bounded  together
with their derivatives. Such solutions may have infinite energy. The proof does not rely on energy estimates, but
on a combination of the maximum principle and quantitative Schauder estimates. We obtain precise bounds on the
sup norm of the solution and its derivatives, making it plain that there is  no exponential increase in time. In particular, these bounds
are time-independent if $g$ is zero. To get a classical solution, it suffices to assume that the initial condition
and the forcing term have bounded derivatives up to order two.
}
\end{quote}

\vspace{4mm}
\noindent

 \medskip
 \noindent {\bf Keywords:}
  viscous Burgers equation, conservation laws,  maximum principle, Schauder estimates.

\smallskip

\noindent
{\bf Mathematics Subject Classification (2010):}   35A01, 35B45, 35B50, 35K15, 35Q30, 35Q35, 35L65, 76N10.

\newpage

\tableofcontents



\section{Introduction and scheme of proof}


\subsection{Introduction}


The $(1+d)$-dimensional viscous Burgers equation is the following non-linear PDE, 
\BEQ (\partial_t -\nu\Del+u\cdot\nabla )u=g, \qquad u\big|_{t=0}=u_0 \label{eq:Burgers}\EEQ
for a velocity $u=u(t,x)\in\R^d$ ($d\ge 1$), $(t,x)\in\R_+\times\R^d$,
where $\nu>0$ is a viscosity coefficient, $\Del$ the standard Laplacian on $\R^d$,   $u\cdot\nabla u=\sum_{i=1}^d
u_i\partial_{x_i}u$ the convection term, and $g$ a  continuous forcing term.  Among other things, this
fluid equation describes the hydrodynamical limit of interacting particle systems \cite{Spo,KipLan}, is a
simplified version without pression of the incompressible Navier-Stokes equation, and also (assuming $g$ to
be random) an interesting toy model for the study of turbulence \cite{BecKha}. The present study is purely mathematical:
we show under  the  following set of   assumptions on
 $u_0$ and $g$ that
 the Cauchy problem 
 \BEQ (\partial_t -\nu\Del+u\cdot\nabla )u=g, u\big|_{t=0}=u_0 \EEQ
  has a unique, globally defined, classical
 solution in $C^{1,2}$ (i.e. continuously differentiable in the time coordinate and twice
 continuously differentiable in the space coordinates), and provide explicit
 bounds for the supremum of $u$ and its derivatives up to second order. 

\medskip

\noindent {\bf Assumptions.} 
{\em \begin{itemize}
\item[(i)] (initial condition) $u_0\in C^{2}$ and $\nabla^2 u_0$ is $\alpha$-H\"older  for every $\alpha\in(0,1)$;   for $\kappa=0,1,2$, $||\nabla^{\kappa}u_0||_{\infty}:=\sup_{x\in\R^d}|\nabla^{\kappa}u_0(x)|<\infty;$
\item[(ii)] (forcing term) on every subset $[0,T]\times\R^d$ with $T>0$ finite, $g$ is bounded and  $\alpha$-H\"older continuous for every $\alpha\in(0,1)$; furthermore, $g$ is $C^{1,2}$ and
  $t\mapsto ||\nabla^{\kappa} g_t||_{\infty}:=\sup_{x\in\R^d} |\nabla^{\kappa} g_t(x)|$,
  $t\mapsto ||\partial_t g_t||_{\infty}:=\sup_{x\in\R^d} |\partial_t g_t(x)|$  are locally
integrable in time. 
\end{itemize}}

For convenience we redefine $\tilde{t}=\nu t$, $\tilde{u}=\nu^{-1} u$, $\tilde{g}=\nu^{-2}g$. The rescaled
equation, $(\partial_{\tilde{t}}-\Del-\tilde{u}\cdot\nabla) \tilde{u}=\tilde{g}$, has viscosity $1$. We skip
the tilde in the sequel. Our bounds blow up in the vanishing viscosity limit
$\nu\to 0$ (see Remarks after Theorem \ref{th:main} for a precise statement).

\bigskip

Our approach is the following. We solve inductively the 
linear transport equations,
\BEQ u^{(-1)}:=0; \EEQ
\BEQ (\partial_t-\Del+u^{(m-1)}\cdot\nabla)u^{(m)}=g, \ \ u^{(m)}\big|_{t=0}=u_0 \qquad (m\ge 0). \label{eq:un}\EEQ
If the sequence $(u^{(m)})_m$ converges in appropriate norms, then the limit is a fixed point of (\ref{eq:un}),
hence solves the Burgers equation. Let $||\ ||_{\alpha}$ denotes either the  isotropic H\"older semi-norm on $\R^d$, $||u_0||_{\alpha}:=\sup_{x,y\in\R^d} \frac{|u_0(x)-u_0(y)|}{|x-y|^{\alpha}}$, or the parabolic H\"older semi-norm on $\R_+\times\R^d$,
$||g||_{\alpha}:=\sup_{(s,x),(t,y)\in\R_+\times\R^d} \frac{|g(s,x)-g(t,y)|}{|x-y|^{\alpha}
+|t-s|^{\alpha/2}}$ (see section 4 for more on H\"older norms).

\begin{Definition}
Let, for $c>0$,  
\BEQ K_0(t):=||u_0||_{\infty}+\int_0^t ds ||g_s||_{\infty} \label{eq:1.5}\EEQ
\BEQ K_1(t):=||\nabla u_0||_{\infty}+\int_0^t ds ||\nabla g_s||_{\infty} \label{eq:1.6}\EEQ
\BEQ  K_2(t):=||\nabla^2 u_0||_{\infty}+||u_0||_{\infty} ||\nabla u_0||_{\infty} + ||g_0||_{\infty}
+\int_0^t ds  \left(||\nabla^2 g_s||_{\infty} + ||\partial_s g_s||_{\infty} \right) \EEQ
\BEQ K_{2+\alpha}(t):=||\nabla^2 u_0||_{\alpha} +  ||g_s||_{\alpha,[0,t]\times\R^d},
\qquad \alpha\in(0,1) \EEQ
and
\BEQ K(t):=c^2 \left(K_0(t)^2+K_1(t)+K_2(t)^{2/3}+K_{2+\alpha}(t)^{2/(3+\alpha)} \right).\label{eq:1.9} \EEQ
\end{Definition}

Note that $K_0(t),K_1(t),K_2(t),K_{2+\alpha}(t),K(t)<\infty$ for all $t\ge 0$ and
$\alpha\in(0,1)$ under the above Assumptions.

\medskip

Our main result is 
the following.

\begin{Theorem} \label{th:main}
For every $\beta\in (0,\half)$, there exists an absolute constant $c=c(d,\beta)\ge 1$, depending only on the dimension and on the exponent $\beta$, such that
the following holds.

\begin{itemize}
\item[(i)] (uniform estimates) 
\BEQ ||u_t^{(m)}||_{\infty}\le K_0(t), \ \  ||\nabla u_t^{(m)}||_{\infty}\le K(t); 
\qquad  ||\partial_t u_t^{(m)}||_{\infty}, ||\nabla^2 u_t^{(m)}||_{\infty} \le (cK(t))^{3/2} \label{eq:uniform-estimates} \EEQ
\item[(ii)] (short-time estimates) define $v^{(m)}:=u^{(m)}-u^{(m-1)}$ for  $m\ge 1$. If $0\le t\le T$ and $t\le m/cK(T)$, then 
\BEQ ||v_t^{(m)}||_{\infty}\le cK_0(T) (cK(T)t/m)^{ m}, \qquad 
||\nabla v_t^{(m)} ||_{\infty}\le cK(T) (cK(T)t/m)^{\beta m}. \label{eq:short-time-estimates}\EEQ
\end{itemize}
\end{Theorem}

Let us comment on these estimates.

\begin{enumerate}
\item The different powers in the expression of $K(t)$ come from the dimension counting dictated by the Burgers equation:
the diffusion term $\Del u$, the convection term $u\cdot\nabla u$ and the forcing $g$ are homogeneous if $u$ scales like $L^{-1}$,
where $L$ is a reference space scale, and $g$ like $(LT)^{-1}$, where $T$ is a reference time scale. Assuming parabolic scaling, $K^{-1}(t)$ scales like time and plays the
r\^ole of a reference time scale $T(t)$ at time $t$, leading to a time-dependent space scale $L=L(t)\sim K^{-\half}(t)$.
The scaling of the other $K$-parameters is $K_0\sim T^{-\half}$; $K_1,K\sim T^{-1}$; $K_2\sim T^{-3/2}$; 
$K_{2+\alpha}\sim T^{-(3+\alpha)/2}$. 

\item The first uniform estimate 
\BEQ ||u_t^{(m)}||_{\infty}\le K_0(t) \EEQ
 follows from a straightforward application of the maximum principle to the transport equation (\ref{eq:un}).
\item (uniform estimates for the gradient).  The function $u^{(0)}$ satisfies the linear heat equation $(\partial_t-\Del)u^{(0)}=g$, whose explicit solution is
$u^{(0)}(t)=e^{t\Del}u_0+\int_0^t ds e^{(t-s)\Del}g_s$. Thus
\BEQ ||\nabla u_t^{(0)}||_{\infty}\le ||\nabla u_0||_{\infty}+\int_0^t ds ||\nabla g_s||_{\infty}=K_1(t). \label{eq:0.11}\EEQ
Clearly $K_1(t)\le K(t)$. Estimates for
further iterates $u^{(1)},u^{(2)},\ldots$ involve $K(t)$ instead of $K_1(t)$.

\item Fix a time horizon $T>0$ and consider the series $S(t):=\sum_{m=0}^{+\infty} v_t^{(m)}=\sum_{m=0}^{+\infty} (u^{(m)}_t-u^{(m-1)}_t)$ for 
$t\le T$ (note that, by  definition, $v^{(0)}:=u^{(0)}-u^{(-1)}=u^{(0)}$). The short-time estimates (\ref{eq:short-time-estimates}) imply that $S(t)$ is absolutely convergent.
More precisely, letting $m_0:=\lfloor cK(T)t\rfloor$ and $\gamma:=1$, 
\BEA ||u_t^{(n)}||_{\infty}=\left|\left| \sum_{m=0}^{n} (u^{(m)}_t-u^{(m-1)}_t) \right|\right|_{\infty} &\le & ||u_t^{(m_0)}||_{\infty} +
\sum_{m=m_0+1}^{+\infty} ||v_t^{(m)}||_{\infty} \nonumber\\
&\le & K_0(T) \left\{ 1+c \sum_{m=m_0+1}^{+\infty} (cK(T)t/m)^{\gamma m} \right\}\EEA
for all $n\ge m_0$.
Let $m>m_0$ and $x=1-cK(T)t/m\in[0,1]$: using $1-x\le e^{-x}$, one gets $(cK(T)t/m)^{\gamma m}=(1-x)^{\gamma m}\le e^{\gamma cK(T)t} e^{-\gamma m}$ and
\BEQ \sum_{m=m_0+1}^{+\infty} (cK(T)t/m)^{\gamma m}\le e^{\gamma cK(T)t} \sum_{m=m_0+1}^{+\infty} e^{-\gamma m}\le e^{\gamma}/(e^{\gamma}-1).\EEQ
 Hence
$||u_t^{(n)}||_{\infty} \lesssim K_0(T).$
In a similar way, letting $\gamma:=\beta$ this time, one shows that
 \BEQ ||\nabla u^{(n)}_t||_{\infty}=\left|\left| \sum_{m=0}^{n} (\nabla u^{(m)}_t-\nabla u^{(m-1)}_t) \right|\right|_{\infty} \lesssim K(T).\EEQ
   These estimates are best when $t=T$; one then retrieves the
   uniform estimates (\ref{eq:uniform-estimates}) up to some constant.

\item (short-time estimates) Bounds (\ref{eq:short-time-estimates}) are of order
$O((Ct)^{\gamma m}/(m!)^{\gamma})$, $\gamma=1$ or $\beta$, and obtained by $m$ successive integrations. For linear
equations, or equations with bounded, uniformly Lipschitz coefficients, successive integrations typically yield $O((Ct)^m/m!)$. The Burgers equation, on the other hand is strongly non-linear. While using precise Schauder estimates
to obtain the gradient bound in  (\ref{eq:short-time-estimates}), one stumbles into the condition $\beta<\half$ at the very
end of section 3 which apparently cannot be improved.

\item (blow-up of the above estimates in the vanishing viscosity limit) Undoing the
initial rescaling, we obtain $\nu$-dependent estimates,
\BEQ ||u_t||_{\infty}\le K_0(t),\qquad ||\nabla u_t||_{\infty}\lesssim \nu^{-1} K(t),
\qquad ||\partial_t u_t||_{\infty}\lesssim \nu^{-1}{K}(t)^{3/2},\qquad 
||\nabla^2 u_t||_{\infty}\lesssim \nu^{-2} {K}(t)^{3/2}\EEQ
with $K_0(t),K_1(t)$ as in (\ref{eq:1.5}), (\ref{eq:1.6}), 
$K_2(t):=\nu||\nabla^2 u_0||_{\infty}+||u_0||_{\infty} ||\nabla u_0||_{\infty}+
||g_0||_{\infty}+\int_0^t ds (\nu||\nabla^2 g_s||_{\infty}+||\partial_s g_s||_{\infty}),$
$K_{2+\alpha}(t):=\nu||\nabla^2 u_0||_{\alpha}+\sup_{\alpha\in[0,t]} ||g_s||_{\alpha}$
and $K(t):=K_0(t)^2+\nu K_1(t)+(\nu K_2(t))^{2/3}+(\nu^{1+\alpha} K_{2+\alpha}(t))^{2/(3+\alpha)}$. Thus the derivative bounds $||\nabla^{\kappa}u_t||_{\infty}$, $\kappa=1,2$
and $||\partial_t u||_{\infty}$ blow up at different rates when $\nu\to 0$.
 
\end{enumerate}

From the above theorem, one deduces easily that the solution of the Burgers equation is smooth on $\R_+\times\R^d$
provided  (i) $u_0$ is smooth and its derivatives are bounded; (ii) $g$ is smooth and its
derivatives  are bounded on $[0,T]\times\R^d$ for all $T$:

\begin{Corollary} \label{cor:main}
Assume $u_0$ and $g$ are smooth, and $||\nabla^{\kappa}u_0||_{\infty}<\infty$ ($\kappa=0,1,2,\ldots$),
$||\partial_t^{\mu}\nabla^{\kappa}g_t||_{\infty}<C(\mu,\kappa,T)$, $\mu,\kappa=0,1,2,\ldots$ for
every $t\le T$. Then the Burgers equation (\ref{eq:Burgers}) has a unique smooth solution $u$ such that
$||\partial_t^{\mu}\nabla^{\kappa}u_t||_{\infty}<C'(\mu,\kappa,T)$ for every $\mu,\kappa$ and $t\le T$. In particular, $C'(\mu,\kappa,t)=C'(\mu,\kappa)$ is uniform in time if
$g=0$.
\end{Corollary}

We do not prove this corollary, since it results from standard extension to higher-order derivatives of the
initial estimates of section 2, and an equally standard iterated use of Schauder estimates to derivatives of
Burgers equation.

\medskip

Our results extend without any modification to nonlinearities of the type $F(u)\cdot\nabla u$ with smooth matrix-valued
coefficient $F$ if $F$ is sublinear, and even (with different scalings and exponents for the $K$-constants) to the
case when $F$ has polynomial growth at infinity. 

\bigskip

Let us compare with the results available in the literature. The one-dimensional case 
$d=1$ or the  irrotational $d$-dimensional case with $g=\nabla f$ of gradient form, is exactly solvable
through the Cole-Hopf transformation $u=\nabla \log\phi$ which reduces it to a scalar, linear PDE 
$\partial_t \phi=\nu\Del\phi+f\phi$; note also that $\log\phi$ is a solution of the
KPZ (Kardar-Parisi-Zhang) equation. In that case the equation is immediately shown
to be well-defined for every $t>0$ under our hypotheses, and estimates similar to
ours are easily obtained; specifically in $d=1$, an invariant measure is known to exist if $g$ is e.g. a space-time white noise \cite{DDT}. For periodic solutions on the
torus in one dimension, the above results extend to the vanishing viscosity limit
\cite{EKMS}. The reader may refer e.g.  to \cite{Der} for a more extended bibliography.

  So our result is mostly interesting  for $d\ge 2$;
as mentioned above, our scheme of proof extends to more general non-linearities of the form $F(u)\cdot\nabla u$, for which
 the equation is not exactly solvable in general. In this setting, the classical
result is that due to Kiselev and Ladyzhenskaja \cite{KL}. The authors consider solutions
in Sobolev spaces and use repeatedly energy estimates. They work on a bounded domain
$\Omega$ with
Dirichlet boundary conditions, but their results extend with minor modifications to the case
$\Omega=\R^d$.  If $u_0\in{\cal H}^s$ with $s>d/2$, then 
$||u_0||_{\infty}<\infty$ by Sobolev's imbedding theorem. Then the maximum principle gives
$||u_t||_{\infty}\le ||u_0||_{\infty}$ as long as the solution is classical; this key
estimate allows one to bootstrap and get bounds for higher-order Sobolev spaces which
increase exponentially in time, e.g. $||u_t||_{H^1}=O(e^{c||u_0||^2_{\infty}t})$,  as follows
 from the proof of Lemma 3 in \cite{KL}. Compared to these estimates, ours present two
essential improvements: (i) we do not assume any decrease of the data at spatial infinity,
so that they do not necessarily belong to Sobolev spaces; (ii) more importantly perhaps,
our bounds do not increase exponentially in time; in the case the right-hand side $g$
vanishes identically, they are even uniform in time, $K_0(t),K(t)\le C$ where $C$
is a constant depending only on the initial condition.


\subsection{Scheme of proof}


Recall that we  solve inductively the following 
linear transport equations, see (\ref{eq:un}),
\BEQ u^{(-1)}:=0; \EEQ
\BEQ (\partial_t-\Del+u^{(m-1)}\cdot\nabla)u^{(m)}=g, \ \ u^{(m)}\big|_{t=0}=u_0 \qquad (m\ge 0). \EEQ
Under the first set of assumptions, standard results on linear equations show that $u^{(m)}$, $m\ge 0$ is
$C^{1,2}$. {\em Assume} we manage to prove locally uniform convergence of $u^{(m)},\nabla u^{(m)}, \nabla^2 u^{(n)}$
when $m\to\infty$. Then
there exists $u\in C^{1,2}$ such that locally uniformly $u^{(m)}\to u$, $\nabla u^{(m)}\to \nabla u$, 
$\nabla^2 u^{(m)}\to \nabla^2 u$ and $\partial_t u^{(m)}\to \partial_t u$. Hence $\partial_t u^{(m)}=\Del u^{(m)}-u^{(m-1)}\cdot\nabla u^{(m)}+g$
converges locally uniformly to  $\Del u-u\cdot\nabla u+g$, and $\partial_t u=\lim_{m\to\infty}
\partial_t u^{(m)}= \Del u-u\cdot\nabla u+g$. In other words, the limit $u$ is a $C^{1,2}$ solution of the
Burgers equation. 

The key point in our scheme is to prove locally uniform convergence of $u^{(m)}$ and $\nabla u^{(m)}$, and to show uniform bounds in H\"older norms for second order
derivatives $\nabla^2 u^{(m)},\ \partial_t u^{(m)}$; a simple argument (see below)
yields then the convergence of second order derivatives, allowing to apply the
above elementary argument. The basic idea is to rewrite $u$ as $\sum_{m=0}^{+\infty} v^{(m)}$, with $v^{(m)}:=u^{(m)}-u^{(m-1)}$,
and to show that the series is convergent, uniformly in space and locally uniformly in time. 

\medskip

In the sequel we fix a constant $c\ge 1$ such that Theorem \ref{th:main} holds and let
\BEQ \bar{K}_0(t):=cK_0(t),\quad \bar{K}_1(t):=cK_1(t),\quad \bar{K}(t):=cK(t) \EEQ
to simplify notations.

\bigskip
\bigskip

The proof relies on two main ingredients: {\em a priori estimates} coming from the
maximum principle; and {\em Schauder estimates}. Schauder estimates are difficult
to find in a precise form suitable for the kind of applications we have in view, so
the reader will find in the  appendix a precise version of these estimates, see
Proposition \ref{prop:XJW}, following a multi-scale proof
introduced by X.-J. Wang. These imply in particular the following.

\begin{Lemma} \label{lem:bound-nabla2um}
 Let $0\le t\le T$. Then
 \BEQ ||\partial_t u^{(m)}||_{\alpha,[0,T]\times\R^d}, ||\nabla^2 u^{(m)}||_{\alpha,[0,T]\times\R^d}\le \bar{K}(T)^{(3+\alpha)/2}.\EEQ
\end{Lemma}

Lemma \ref{lem:bound-nabla2um} is proved in section 3, at the same time as  Theorem \ref{th:main}.

\medskip

We now use a classical result about H\"older spaces: let $C^{\alpha}(Q)$, with $Q\subset\R\times
\R^d$ compact, be the
Banach space of $\alpha$-H\"older functions on $Q$ equipped with the norm
$|||u|||_{\alpha}:=||u||_{\infty,Q}+||u||_{\alpha,Q}$. Then the injection $C^{\alpha'}(Q)\subset C^{\alpha}(Q)$ is compact for every $\alpha'<\alpha$. In particular, Lemma \ref{lem:bound-nabla2um} implies the existence
of a subsequence $(u^{(n_m)})_m$ such that $\nabla^2 u^{(n_m)}\to_{m\to\infty} v$
in $C^{\alpha'}$-norm. On the other hand, as discussed in Remark 4 above, 
$u^{(m)}\to u$ and $\nabla u^{(m)}\to \nabla u$ in the sup norm for some $u\in C^{0,1}$. Hence $u$ is twice continuously differentiable in the space variables, and $\nabla^2 u=v$. Now every subsequence $(\nabla^2 u^{(n'_m)})_m$ converges
to the same limit, $\nabla^2 u$. Hence $\nabla^2 u^{(n)}\to \nabla^2 u$ in $C^{\alpha'}$.
In a similar way, one proves that $u$ is continuously differentiable in the time
variable, and $\partial_t u=\lim_{m\to\infty} \partial_t u^{(m)}$ in $C^{\alpha'}$.
In particular, $u\in C^{1,2} $, and the arguments given at the very beginning of the present
subsection show that $u$ is a classical solution of the Burgers equation.
Note that we may reach the same conclusion even if we do not know that the series
$||\nabla u^{(m+1)}-\nabla u^{(m)}||_{\infty,Q}$ converges. Actually the bound on
$||\nabla u^{(m+1)}-\nabla u^{(m)}||_{\infty,Q}$ is the trickiest one. We felt however
it was one the most inexpected estimates we had obtained, and thus worth including.

\bigskip

{\em Notations.} For $f,g:X\to\R_+$ two positive functions on a set $X$,  we write $f(u)\lesssim g(u)$ if there exists a constant $C=C(d)$ depending only on
the dimension such that $f(u)\le Cg(u)$. (If $C$ depends on other parameters, notably on $c$, then we write explicitly the 
dependence on them, so that we make it clear that we do not get unwanted extra multiplicative factors $O(c^m)$ in
the formulas which would invalidate the proofs).  


\section{Initial estimates}


Initial estimates are different in spirit from those of the next section since they
cannot rely on Schauder estimates. Instead we use a Gronwall-type lemma based
on the maximum principle.

\begin{Lemma}[Gronwall lemma] \label{lem:Gronwall}
Let $\phi:\R_+\times\R^d\to\R^d$, resp. $\bar{\phi}:\R_+\times\R^d\to\R^d$ be the solution of the transport equation
$(\partial_t-\Del+b\cdot\nabla-c)\phi=f$, resp. $(\partial_t-\Del+\bar{b}\cdot\nabla-\bar{c})\bar{\phi}=\bar{f}$, with same initial
condition, $\phi\big|_{t=0}=\bar{\phi}\big|_{t=0}$; the coefficients $c=c(t,x),\bar{c}=\bar{c}(t,x)\in M_{d\times d}(\R)$ are
 matrix-valued, and $b,\bar{b},c,\bar{c}$ are assumed to be bounded and continuous.  Let $v:=\bar{\phi}-\phi$. Then
\BEQ ||v_t||_{\infty}\le \int_0^t ds\, A(s,t)\,  ||\bar{b}_s-b_s||_{\infty} \,
||\nabla \phi_s||_{\infty} + \int_0^t ds\,  A(s,t)\,  |||\bar{c}_s-c_s|||_{\infty} \, ||\phi_s||_{\infty}
+ \int_0^t ds\, A(s,t)\,  ||\bar{f}_s-f_s||_{\infty},\EEQ
where  $|||\ \cdot \ |||_{\infty}$ is the supremum over $\R^d$ of the operator norm in $M_{d\times d}(\R)$, and 
$A(s,t)=\exp \int_s^t |||\bar{c}_r|||_{\infty} dr.$
\end{Lemma}

{\bf Proof.} By subtracting the PDEs satisfied by $\phi$ and $\bar{\phi}$, one gets
\BEQ (\partial_t-\Del+\bar{b}\cdot\nabla-\bar{c})v=-(\bar{b}-b)\cdot \nabla\phi
+(\bar{f}-f)+(\bar{c}-c)\phi. \label{eq:2.2} \EEQ Hence the result by the maximum principle. \hfill\eop

\begin{Definition} \label{def:tinit}
Let $t_{init}:=\inf\left\{t>0;\ t\bar{K}(t)=1\right\}.$
\end{Definition}

By hypothesis, $t_{init}>0$. If
$u_0\equiv 0$ and $g\equiv 0$, then 
$t_{init}=+\infty$ and the solution of Burgers' equation is simply $0$. The case
$u_0$=Cst, $\nabla g=0$ reduces to the previous one by the generalized Galilean
transformation $x\mapsto x+\int_0^t a(s)ds$, $u\mapsto u-a$ with $a(t)=u_0+\int_0^t g_s ds$. We
henceforth exclude this trivial case, so that $t_{init}\in(0,+\infty)$.

\begin{Theorem}[initial estimates] \label{th:initial}
Let $t\le t_{init}$. Then the following estimates hold:
\begin{itemize}
\item[(i)] 
\BEQ ||u_t^{(m)}||_{\infty}\le K_0(t_{init}), \quad ||\nabla u_t^{(m)}||_{\infty}
\le K(t_{init}); \qquad ||\partial_t u_t^{(m)}||_{\infty}, ||\nabla^2 u_t^{(m)}||_{\infty}\le \bar{K}(t_{init})^{3/2}.
\label{eq:initial-estimates1} \EEQ
Furthermore, 
\BEQ ||\partial_t u_t^{(m)}||_{\alpha}, ||\nabla^2 u_t^{(m)}||_{\alpha} \le C\bar{K}(t_{init})^{(3+\alpha)/2} \label{eq:initial-alpha} \EEQ
with $C=C(d,\alpha)$.
\item[(ii)]  let $m\ge 1$,  then 
\BEQ ||v_t^{(m)}||_{\infty}\le \bar{K}_0(t_{init}) (\bar{K}(t_{init})t/m)^m, \qquad 
||\nabla v_t^{(m)} ||_{\infty}\le \bar{K}(t_{init}) (\bar{K}(t_{init})t/m)^m. \label{eq:initial-estimates2} \EEQ
\end{itemize}
\end{Theorem}

{\bf Remarks.}
\begin{enumerate}
\item Let $T\le t_{init}$, then (\ref{eq:initial-estimates1}), (\ref{eq:initial-alpha}) and (\ref{eq:initial-estimates2}) remain true
for $t\le T$ if one replaces $K_0(t_{init})$, $\bar{K}_0(t_{init})$, $K(t_{init})$, $\bar{K}(t_{init})$ by
$K_0(T)$, $\bar{K}_0(T)$, $K(T),\bar{K}(T)$. Hence Theorem \ref{th:main} is proved for $t\le t_{init}$ (actually
with $\beta=1$).
\item The value of $t_{init}$ depends on the choice of $c$. We provide in the course of
the proof a rather explicit minimal value of $c$ for which (\ref{eq:initial-estimates1}), (\ref{eq:initial-alpha}), (\ref{eq:initial-estimates2}) hold. Further estimates
in the next section may require a larger value of $c$.
\item From H\"older  interpolation estimates (see Lemma \ref{lem:interpolation-Holder}), one also has a bound for lower-order H\"older norms,
\BEQ ||u^{(m)}||_{\alpha}\lesssim  K_0(t_{init})^{1-\alpha}\bar{K}(t_{init})^{\alpha}+
K_0^{1-\alpha/2}(t_{init})\bar{K}^{3\alpha/4}(t_{init}),\EEQ
and, for fixed $s\le t_{init}$,
\BEQ ||\nabla u_s^{(m)}||_{\alpha}\lesssim  K^{1-\alpha}(t_{init})\bar{K}(t_{init})^{3\alpha/2}.\EEQ 

\end{enumerate}

\medskip

{\bf Proof.} Let us abbreviate $K_0(t_{init}),\bar{K}_0(t_{init}), K_1(t_{init}), \bar{K}_1(t_{init})$, $K(t_{init})$, $\bar{K}(t_{init})$ to
$K_0,\bar{K}_0,K_1,\bar{K}_1,K,\bar{K}$.

\begin{itemize}
\item[(i)] We first prove estimates (i)  by induction, assuming them to be
proved for $m-1$. Note first that (\ref{eq:initial-estimates1})  holds true for $m=0$ with $c=1$, see eq. (\ref{eq:0.11}); as for (\ref{eq:initial-alpha}), 
\BEA  ||\nabla^2 u_t^{(0)}||_{\gamma} &\lesssim&  ||\nabla^2 u_0||_{\gamma}+ \int_0^t ds \, ||\nabla^2 e^{s\Del} g_{t-s}||_{\gamma} \nonumber\\
& \le& K_2^{1-\gamma/\alpha}(t_{init}) K_{2+\alpha}^{\gamma/\alpha}(t_{init}) + t_{init}^{(\alpha-\gamma)/2} K_{2+\alpha} (t_{init})\nonumber\\
&\le &   C(d,\alpha,\gamma) \bar{K}^{(3+\gamma)/2}, \qquad \gamma<\alpha \EEA
as follows from H\"older interpolation inequalities (see Lemma \ref{lem:interpolation-Holder}) and Corollary \ref{cor:Holder}. Time variations of $\nabla^2 u_t^{(0)}$ scale
similarly, yielding $||\nabla^2 u^{(0)}||_{\gamma,[0,t_{init}]\times\R^d} \lesssim
\bar{K}^{(3+\gamma)/2}$ (see Lemma \ref{lem:Holder}, eq. (\ref{eq:bound-nabla2t}), and
Corollary \ref{cor:Holder}). Note that similarly, $||\nabla u^{(0)}||_{\gamma,[0,t_{init}]\times\R^d}\lesssim \bar{K}^{(2+\gamma)/2}$. The estimate for $||u_t^{(m)}||_{\infty}$
is a direct consequence of the maximum principle. Then $\nabla u^{(m)}$ satisfies the gradient
equation 
\BEQ (\partial_t-\Del+u^{(m-1)}\cdot\nabla +\nabla u^{(m-1)})\nabla u^{(m)}=
\nabla g, \EEQ
 where $\nabla u^{(m-1)}(t,x)$ is viewed as the $d\times d$ matrix $(\partial_j u_k(t,x))_{jk}$ acting on the vector $(\partial_k u_i)_k$. Note that
\BEQ |||\nabla u^{(m-1)}(t,x)|||\le \sqrt{\Tr(\nabla u^{(m-1)}(t,x))(\nabla u^{(m-1)}(t,x))^* }=|\nabla u^{(m-1)}(t,x)|.
\EEQ
  By the maximum principle,
\BEQ ||\nabla u_t^{(m)}||_{\infty}\le A(0,t)\,  ||\nabla u_0||_{\infty} + \int_0^t ds\, 
A(s,t) \,  ||\nabla g_s||_{\infty},\EEQ
where $A(s,t):= \exp \int_s^t ||\nabla u^{(m-1)}_r||_{\infty} dr$ is the exponential
amplification factor of Lemma \ref{lem:Gronwall}. By induction hypothesis and Definition \ref{def:tinit},
 $A(s,t)\le A(0,t_{init})
\le e^{t_{init}K}\le e$, hence (provided $c^2\ge e$)
\BEQ ||\nabla u_t^{(m)}||_{\infty}\le e K_1\le K.\EEQ

\medskip

To bound $\nabla^2 u_t^{(m)}$ we differentiate once more,
\BEQ (\partial_t-\Del+u^{(m-1)}\cdot\nabla +\nabla u^{(m-1)})\nabla^2 u^{(m)}=
\nabla^2 g- \nabla^2 u^{(m-1)} \nabla u^{(m)},\EEQ
where $\nabla u^{(m-1)}$ is viewed this time as the $d^2\times d^2$ matrix $\left( \partial_{j'}u_k^{(m-1)}\del_{k',j}+\partial_j u_k^{(m-1)}\del_{k',j'}\right)_{(jj'),(kk')}$
acting on the vector  $(\partial^2_{kk'}u_i)_{kk'}\in\R^{d^2}$, and has matrix norm $|||\nabla u^{(m-1)}(t,x)|||_{M_{d^2\times d^2}(\R)}\le C_d |\nabla u^{(m-1)}(t,x)|$, 
yielding an amplification factor $\tilde{A}(s,t):=\exp \int_s^t ||\ |||\nabla u_r^{(m-1)}(t,x)|||_{M_{d^2\times d^2}(\R)} \ ||_{\infty} dr\le C'_d$.
 By the maximum principle,
\BEA   ||\nabla^2 u_t^{(m)}||_{\infty} &\le & C'_d \left( ||\nabla^2 u_0||_{\infty} + \int_0^t ds
 \left( ||\nabla^2 g_s||_{\infty} + ||\nabla^2 u_s^{(m-1)}||_{\infty}
||\nabla u_s^{(m)}||_{\infty} \right) \right) \nonumber\\
&\le & C'_d\left(||\nabla^2 u_0||_{\infty}+\int_0^t ds ||\nabla^2 g_s||_{\infty}+ t_{init}\bar{K}^{3/2}K \right) \nonumber\\
&\le & C'_d(K_2(t_{init})+\bar{K}^{\half}K)\le C'_d(c^{-3}+c^{-1})\bar{K}^{3/2} \le \bar{K}^{3/2}\EEA
provided $c\ge 2 \max(1,C'_d)$.

\medskip
Similarly, $\partial_t u^{(m)}$ satisfies the transport equation
\BEQ (\partial_t-\Del+u^{(m-1)}\cdot \nabla)\partial_t u^{(m)}=\partial_t g- \partial_t u^{(m-1)}\cdot \nabla u^{(m)},
\EEQ
hence
\BEA ||\partial_t u_t^{(m)}||_{\infty} & \le & ||\nabla^2 u_0||_{\infty}+||u_0||_{\infty} \, ||\nabla u_0||_{\infty}
+||g_0||_{\infty}  +\int_0^t ds ||\partial_s g_s||_{\infty}+t_{init}\bar{K}^{3/2}K \nonumber\\
 & \le &  K_2(t_{init})+\bar{K}^{\half}K\le  (c^{-3}+c^{-1})\bar{K}^{3/2}
 \le \bar{K}^{3/2} \EEA
 provided $c\ge 2$.
 
 \medskip
 Finally, we must prove the H\"older estimate (\ref{eq:initial-alpha}): for that, we use the integral representation
 \BEQ \nabla^2 u_t^{(m)}=\nabla^2 u_t^{(0)}-\int_0^t \nabla^2 e^{(t-s)\Del}\left( (u_s^{(m-1)}\cdot\nabla)u_s^{(m)}
 \right)ds.\EEQ
 By Lemma \ref{lem:interpolation-Holder}, considering $\alpha$-H\"older norms on
 $[0,t_{init}]\times\R^d$,
 \BEA  ||(u_s^{(m-1)}\cdot\nabla)u_s^{(m)}||_{\gamma} &\le & ||u_s^{(m-1)}||_{\infty} \ ||\nabla u_s^{(m)}||_{\gamma}
 + ||\nabla u_s^{(m)}||_{\infty} \ ||u_s^{(m-1)}||_{\gamma} \nonumber\\
 &\lesssim& K_0 K^{1-\gamma} \bar{K}^{3\gamma/2} + K K_0^{1-\gamma}\bar{K}^{\gamma} \lesssim  \bar{K}^{(3+\gamma)/2} \EEA
 Thus by Lemma \ref{lem:Holder},
 \BEA  ||\nabla^2 u_t^{(m)}-\nabla^2 u_{t'}^{(m)}||_{\infty} &  \lesssim &  ||\nabla^2 u_t^{(0)}-\nabla^2 u_{t'}^{(0)}||_{\infty}+ 
 \int_{t'}^t (t-s)^{\frac{\alpha}{2}-1} ||(u_s^{(m-1)}\cdot\nabla)u_s^{(m)}||_{\alpha}ds
\nonumber\\  &\lesssim&  (t-t')^{\alpha/2} \bar{K}^{(3+\alpha)/2} \EEA
 for $t'<t$, and (choosing any $\gamma\in(\alpha,1)$)
 \BEQ ||\nabla^2 u_t^{(m)}||_{\alpha}\lesssim ||\nabla^2 u_t^{(0)}||_{\alpha}+
 C'(d,\alpha,\gamma) \bar{K}^{(3+\gamma)/2} \int_0^{t_{init}} (t-s)^{-1+(\gamma-\alpha)/2}
 ds\lesssim \bar{K}^{(3+\alpha)/2},\EEQ
 hence the result for $||\nabla^2 u^{(m)}||_{\alpha}$.  Similarly, 
 \BEA  ||\nabla u_t^{(m)}-\nabla u_{t'}^{(m)}||_{\alpha} & \lesssim & ||\nabla u_t^{(0)}-
 \nabla u_{t'}^{(0)}||_{\alpha}+ 
 \int_{t'}^t (t-s)^{(\alpha-1)/2} ||(u_s^{(m-1)}\cdot\nabla)u_s^{(m)}||_{\alpha} ds
\nonumber\\ &\lesssim & (t-t')^{\alpha/2} \bar{K}^{(2+\alpha)/2} + (t-t')^{(\alpha+1)/2} \bar{K}^{(3+\alpha)/2}  \nonumber\\ & \lesssim &   (t-t')^{\alpha/2} \bar{K}^{(2+\alpha)/2}+ (t-t')^{\alpha/2} t_{init}^{\half}
 \bar{K}^{(3+\alpha)/2}\lesssim (t-t')^{\alpha/2} \bar{K}^{(2+\alpha)/2}, \nonumber\\ \EEA
 hence (using H\"older interpolation inequalities once more) $||\nabla u^{(m)}||_{\alpha}\lesssim \bar{K}^{(2+\alpha)/2}$. From the previous bounds follows immediately
 $||\partial_t u^{(m)}||_{\alpha}\lesssim ||\nabla^2 u^{(m)}||_{\alpha}+||(u^{(m-1)}\cdot \nabla)u^{(m)}||_{\alpha}\lesssim \bar{K}^{(3+\alpha)/2}$.

\item[(ii)] 
Apply Lemma \ref{lem:Gronwall} with $\phi=\bar{b}=u^{(m-1)}$, $b=u^{(m-2)}$, $\bar{\phi}=u^{(m)}$, $f=\bar{f}=g$ and $c=\bar{c}=0$. It comes out
\BEQ ||v_t^{(m)}||_{\infty}\le \int_0^t ds ||v_s^{(m-1)}||_{\infty} ||\nabla u_s^{(m-1)}||_{\infty}.\EEQ
Thus, using the induction hypothesis,
\BEQ ||v_t^{(m)}||_{\infty}\le \int_0^t ds \bar{K}_0 (\bar{K}s/(m-1))^{m-1} K
\le \bar{K}_0 (\bar{K}t/m)^m (1-\frac{1}{m})^{-(m-1)} (K /\bar{K})\le   \bar{K}_0 (\bar{K}t/m)^m, \qquad m\ge 2 \EEQ
for $c$ large enough, and
\BEQ ||v_t^{(1)}||_{\infty}\le \int_0^t ds ||u_s^{(0)}||_{\infty} ||\nabla u_s^{(0)}||_{\infty} \le K_0 K t \le \bar{K}_0 (\bar{K}t).\EEQ

\medskip

Consider now as in (i) the gradient of the transport equations of index $m-1,m$,
\BEQ (\partial_t-\Del+u^{(n-1)}\cdot\nabla+\nabla u^{(n-1)})\nabla u^{(n)}=\nabla g,
\qquad n=m-1,m \EEQ
and apply  Lemma \ref{lem:Gronwall} with $\phi=\nabla u^{(m-1)}$, $\bar{\phi}=\nabla u^{(m)}$,  $b=u^{(m-2)}$, $\bar{b}=u^{(m-1)}$ and $c=\nabla u^{(m-2)}$, $\bar{c}=\nabla u^{(m-1)}$. Using the induction hypothesis, one gets
\BEA ||\nabla v_t^{(m)}||_{\infty} &\le & \int_0^t ds\,  A(s,t)\,  ||v_s^{(m-1)}||_{\infty} ||\nabla^2 u_s^{(m-1)}||_{\infty} +
 \int_0^t ds\,  A(s,t)\, ||\nabla v_s^{(m-1)}||_{\infty} ||\nabla u_s^{(m-1)}||_{\infty} \nonumber\\
&\le & e \int_0^t ds (\bar{K}_0\bar{K}^{3/2}+\bar{K}K) (\bar{K}s/(m-1))^{m-1}   \nonumber\\
&\le & e (1-\frac{1}{m})^{-(m-1)}  (\bar{K}t/m)^m (\bar{K}_0 \bar{K}^{\half} +K)\le e (1-\frac{1}{m})^{-(m-1)} (c^{-\half}+c^{-1}) \bar{K} (\bar{K}t/m)^m  \nonumber\\
&\le & 
 \bar{K} (\bar{K}t/m)^m, \qquad m\ge 2  \EEA
 and
 \BEA ||\nabla v_t^{(1)}||_{\infty} &\le & \int_0^t ds \left( ||u_s^{(0)}||_{\infty}
 ||\nabla^2 u_s^{(0)}||_{\infty} + ||\nabla u_s^{(0)}||_{\infty}^2 \right) \nonumber\\
 &\le & e(K_0 \bar{K}^{3/2} + K^2)t \le \bar{K} (\bar{K}t) \EEA
for $c$ large enough.

\end{itemize}
\hfill \eop


\section{Proof of main theorem}


By Remark 1  following Theorem \ref{th:initial},   we may now restrict to times larger than $t_{init}$.
We  fix a time horizon $T>t_{init}$ and distinguish two regimes: a 
{\em short-time regime}, $t \le m/\bar{K}(T)$; and a {\em long-time regime}, $t> m/\bar{K}(T)$. Clearly the short-time regime does not exist for $m=0$; 
as already noted before (see comments after Theorem \ref{th:main}), this  case is trivial
and 
estimates (\ref{eq:uniform-estimates}), proven in the course of
Theorem \ref{th:initial} in the initial regime, extend without any modification to
arbitrary time. So we assume henceforth that $m\ge 1$.

\medskip

Theorem \ref{th:main} follows immediately from an estimate for $u^{(m)},\nabla u^{(m)}$ valid
over the whole region $t\in[t_{init},T]$ and another estimate for $v^{(m)},\nabla v^{(m)}$ valid only in the short-time regime. These are proved by induction.

\begin{Theorem}[estimates for $u^{(m)}$ and $\nabla u^{(m)}$] \label{th:u-estimates}
Let $m\ge 1$ and $t \in [t_{init},T]$. Then 
\BEQ ||u_t^{(m)}||_{\infty}\le K_0(T), \ ||\nabla u_t^{(m)}||_{\infty}\le K(T); \qquad
||\partial_t u_t^{(m)}||_{\infty},||\nabla^2 u_t^{(m)}||_{\infty}\le \bar{K}(T)^{3/2}.
\label{eq:u-estimates} \EEQ
Furthermore, 
\BEQ ||\partial_t u_t^{(m)}||_{\alpha}, ||\nabla^2 u_t^{(m)}||_{\alpha}\lesssim
\bar{K}(T)^{(3+\alpha)/2}.\EEQ
\end{Theorem}

{\bf Proof.}  As already noted, the inequality $||u_t^{(m)}||_{\infty}\le K_0(T)$ follows immediately from the maximum principle, so we consider only the
bound for the gradient and higher-order derivatives in (\ref{eq:u-estimates}). We prove it by induction on $m$, assuming it to be true for $m-1$.
We abbreviate $K_0(T),K(T),\bar{K}(T)$ to $K_0,K,\bar{K}$.

\medskip

We apply Proposition \ref{prop:XJW} on the parabolic ball $Q^{(j)}=[t-M^j,t]\times \bar{B}(x,M^{j/2})$,
with $M^j:=\half \bar{K}(T)^{-1}$. Note that, by definition, $t-M^j\ge t_{init}-\half \bar{K}(t_{init})^{-1}\ge
\half t_{init}>0$. We consider first the bound (\ref{eq:3.46}) for the gradient,
\BEA && ||\nabla u^{(m)}||_{\infty,Q^{(j-1)}} \lesssim R_b^{-1} \bar{K}^{-(\alpha+1)/2} ||g||_{\alpha,Q^{(j)}}
+ R_b^{-1} K_0 \left( \bar{K}^{-(\alpha+\half)} R_b^{-1} ||u^{(m-1)}||^2_{\alpha,Q^{(j)}} + \bar{K}^{\half} \right).
\label{eq:uv} \EEA
The multiplicative factor $R_b^{-1}$ is bounded by $1+ (2\bar{K})^{-\half} ||u^{(m-1)}||_{\infty,Q^{(j)}}\le 1+\bar{K}^{-\half} K_0\le 2$. On the other hand, by H\"older interpolation inequalities (see Lemma \ref{lem:interpolation-Holder}), 
\BEA  ||u^{(m-1)}||_{\alpha,Q^{(j)}} &\lesssim&  K^{\alpha} K_0^{1-\alpha}+ \bar{K}^{3\alpha/4} K_0^{1-\alpha/2} \nonumber\\
&\le & (1+c^{3\alpha/4}(K_0^2/K)^{\alpha/4})K^{\alpha}K_0^{1-\alpha} \nonumber\\
&\le & (1+c^{\alpha/4}) K^{\alpha}K_0^{1-\alpha} \le (1+c^{\alpha/4})c^{2\alpha-2} K^{(1+\alpha)/2}. \label{eq:3.4} \EEA 
Hence 
\BEA  ||\nabla u^{(m)}||_{\infty,Q^{(j-1)}} &\lesssim &
\bar{K}^{-\alpha-1/2} K_{2+\alpha}(T) +  K_0 \bar{K}^{-\alpha-\half} \cdot
c^{\alpha/2} K^{2\alpha} K_0^{2-2\alpha} + \bar{K}^{\half} K_0
\nonumber\\
&\le & c^{-\alpha-1/2} K+ c^{-(1+\alpha)/2} K^{\alpha-\half} K_0^{3-2\alpha} + c^{-\half} K  \EEA
which is $\le K$ for $c$ large enough.

\medskip
Bounds for higher-order derivatives $||\partial_t u_t^{(m)}||_{\infty}, ||\nabla^2 u_t^{(m)}||_{\infty}$ follow
from (\ref{eq:3.48}) instead, contributing an extra $M^{-j/2}\approx \bar{K}^{\half}$ multiplicative factor.
They hold true for $c$ large enough. Finally, (\ref{eq:3.49}) yields
\BEA ||\partial_t u^{(m)}||_{\alpha,Q^{(j-1)}},||\nabla^2 u^{(m)}||_{\alpha,Q^{(j-1)}} &\lesssim &
 ||g||_{\alpha,Q^{(j)}} + K_0 \left( ||u^{(m-1)}||_{\alpha,Q^{(j)}}^{(2+\alpha)/(1+\alpha)}
 + \bar{K}^{1+\alpha/2} \right) \nonumber\\
 &\lesssim & K_{2+\alpha}(T)+ K_0 \cdot c^{(\frac{\alpha}{4}+2\alpha-2)(2+\alpha)/(1+\alpha)} K^{1+\alpha/2}
 +c^{-1} \bar{K}^{(3+\alpha)/2}  \nonumber\\
 &\lesssim & \bar{K}^{(3+\alpha)/2}, \EEA

from which 
\BEQ ||\nabla^2 u^{(m)}||_{\alpha,[t_{init},T]\times\R^d} \lesssim \sup_{(t,x)\in[t_{init},T]\times\R^d} ||\nabla^2 u^{(m)}||_{\alpha, Q^{(j-1)}(t,x)}+M^{-j\alpha/2}
||\nabla^2 u^{(m)}||_{\infty,[t_{init},T]\times\R^d}\lesssim \bar{K}^{(3+\alpha)/2},\EEQ
and similarly for $||\partial_t u^{(m)}||_{\alpha,[t_{init},T]\times\R^d}$. 

\medskip
We take the opportunity to derive from (\ref{eq:3.47}) a bound for $||\nabla u^{(m)}||_{\alpha,Q^{(j-1)}}$ (also valid for $||\nabla u^{(m)}||_{\alpha,[t_{init},T]\times\R^d}$) 
that will be helpful in the next theorem,
\BEA &&  ||\nabla u^{(m)}||_{\alpha,Q^{(j-1)}}  \lesssim \bar{K}^{-1/2} (1+\bar{K}^{-(1+\alpha)/2} ||u^{(m-1)}||_{\alpha,
Q^{(j)}}) ||g||_{\alpha} + \nonumber\\
&&\qquad \qquad \qquad \qquad  K_0 \bar{K}^{(1+\alpha)/2} \left( 1+\bar{K}^{-(1+\alpha)/2} ||u^{(m-1)}||_{\alpha,Q^{(j)}}
+ (\bar{K}^{-(1+\alpha)/2} ||u^{(m-1)}||_{\alpha,Q^{(j)}})^3 \right) \nonumber\\
&&\qquad \qquad \qquad \lesssim \bar{K}^{1+\alpha/2} \EEA
since (from (\ref{eq:3.4})) $||u^{(m-1)}||_{\alpha,Q^{(j)}}\lesssim \bar{K}^{(1+\alpha)/2}$.

 \hfill \eop

\begin{Theorem}[short-time estimates for $v^{(m)}$ and $\nabla v^{(m)}$] \label{th:v-estimates}
Let $m\ge 1$ and $t\in [t_{init},\min(T,m/\bar{K}(T))]$. Then 
\BEQ ||v_t^{(m)}||_{\infty}\le \bar{K}_0(T) (\bar{K}(T)t/m)^{ m}, \qquad 
||\nabla v_t^{(m)} ||_{\infty}\le \bar{K}(T) (\bar{K}(T)t/m)^{\beta m}. \label{eq:v-estimates}\EEQ
\end{Theorem}

{\bf Proof.} We abbreviate as before $K_0(T),\bar{K}_0(T)$, $K(T),\bar{K}(T)$ to $K_0$, $\bar{K}_0$,
$K$, $\bar{K}$ and prove simultaneously the bounds on $||v^{(m)}||_{\infty}$ and $||\nabla v^{(m)}||_{\infty}$,
assuming them to be true for $m-1$. 

\begin{itemize}
\item[(i)] (bound for $v_t^{(m)}$) 
As in the proof of Theorem \ref{th:initial} (ii), the case $m=1$ is 
essentially trivial: namely, using Lemma \ref{lem:Gronwall}, we have for $t\le \bar{K}^{-1}$
\BEQ ||v_t^{(1)}||_{\infty}\le \int_0^t ds  \, ||u_s^{(0)}||_{\infty} \, ||\nabla u_s^{(0)}||_{\infty}\le 
  K_0 Kt\le \bar{K}_0 (\bar{K}t).\EEQ So we now restrict to $m\ge 2$.

Assume first $t\le (m-1)/\bar{K}$, so that $t$ is in the short-time regime for $u^{(m-1)}$. By Lemma \ref{lem:Gronwall}
(see proof of Theorem \ref{th:initial} (ii)),
\BEA   ||v_t^{(m)}||_{\infty} & \le & \int_0^t ds  \, ||v_s^{(m-1)}||_{\infty} \, ||\nabla u_s^{(m-1)}||_{\infty} \nonumber\\
&\le & \int_0^t ds  \, \bar{K}_0 (\bar{K}s/(m-1))^{m-1} K\le (\bar{K}t/(m-1))^{m}  \bar{K}_0 (K/\bar{K}) \nonumber\\
&\le & c^{-1}\bar{K}_0 (\bar{K}t/(m-1))^{m}  \le \half \bar{K}_0 (\bar{K}t/m)^{m}   \label{eq:3.8}\EEA
for $c$ large enough. 

\medskip
For $s,t\in[(m-1)/\bar{K},m/\bar{K}]$, one uses instead $||v_s^{(m-1)}||_{\infty}\le ||u_s^{(m-1)}||_{\infty}+
||u_s^{(m-2)}||_{\infty}\le 2K_0$ and obtains
\BEA  ||v_t^{(m)}||_{\infty} & \le & \int_0^{(m-1)/\bar{K}} ds  \, ||v_s^{(m-1)}||_{\infty} \, ||\nabla u_s^{(m-1)}||_{\infty}
+ \int_{(m-1)/\bar{K}}^{m/\bar{K}} ds \, ||v_s^{(m-1)}||_{\infty} \, ||\nabla u_s^{(m-1)}||_{\infty} \nonumber\\
&\le & \half \bar{K}_0 (\bar{K}t/m)^{m} +  \bar{K}^{-1}\cdot 2K_0 K \nonumber\\
&\le & \bar{K}_0 (\bar{K}t/m)^{m}   \label{eq:3.9} \EEA
for $c$ large enough. 

\item[(ii)] (bound for $\nabla v_t^{(m)}$)
We start from the observation (see (\ref{eq:2.2}))  that $v^{(m)}$ satisfies the transport equation 
$(\partial_t-\Del+u^{(m-1)}\cdot\nabla)(v^{(m)})=-v^{(m-1)}\cdot \nabla u^{(m-1)}$ and apply
Schauder estimates on $Q^{(j)}=Q^{(j)}(t_0,x_0)$  as in the proof of Theorem \ref{th:u-estimates}, with $M^j\approx \bar{K}(T)^{-1}$,
$b=u^{(m-1)}$ and $f:=v^{(m-1)}\cdot \nabla u^{(m-1)}$.  In the course of the proof of Theorem \ref{th:u-estimates},
and in (i),
we obtained  $||u^{(m-1)}||_{\infty,Q^{(j)}}\le K_0$ and 
\BEQ   ||u^{(m-1)}||_{\alpha,Q^{(j)}}\lesssim \bar{K}^{(1+\alpha)/2}, 
\qquad ||v^{(m)}||_{\infty,Q^{(j)}}\le \bar{K}_0 (\bar{K}t/m)^m, \qquad ||\nabla u^{(m-1)}||_{\alpha,
Q^{(j)}}\lesssim \bar{K}^{1+\alpha/2}.\EEQ
Furthermore, from H\"older interpolation inequalities (see Lemma \ref{lem:interpolation-Holder}) and induction hypothesis,
\BEQ ||v^{(m-1)}||_{\alpha,Q^{(j)}}\lesssim \bar{K}_0^{1-\alpha} \bar{K}^{\alpha} (\bar{K}t/(m-1))^{\beta(m-1)}.\EEQ
Hence (using once again the induction hypothesis)
\BEA ||f||_{\alpha,Q^{(j)}} &\lesssim &  ||v^{(m-1)}||_{\alpha,Q^{(j)}} \, ||\nabla u^{(m-1)}||_{\infty,Q^{(j)}} +
||v^{(m-1)}||_{\infty,Q^{(j)}} \, ||\nabla u^{(m-1)}||_{\alpha,Q^{(j)}} \nonumber\\
&\lesssim& (\bar{K}t/(m-1))^{\beta(m-1)} (\bar{K}_0^{1-\alpha} \bar{K}^{\alpha}K+\bar{K}_0 \bar{K}^{1+\alpha/2})
 \nonumber\\
 &\lesssim & c^{-1} \bar{K}^{(3+\alpha)/2} (\bar{K}t/(m-1))^{\beta(m-1)}.  \label{eq:falpha} \EEA
 
\medskip

A priori we should now use the Schauder estimate (\ref{eq:3.47}) to bound $||\nabla v^{(m)}||_{\alpha,Q^{(j-1)}}$; as
in the proof of Theorem \ref{th:u-estimates}, $R_b^{-1}\le 2$, so
\BEA  ||\nabla v^{(m)}||_{\infty,Q^{(j-1)}} & \lesssim & \bar{K}^{-(1+\alpha)/2} ||f||_{\alpha} + 
\bar{K}^{1/2} \bar{K}_0 \left(1+ (\bar{K}^{-1-\alpha/2} ||u^{(m-1)}||_{\alpha})^2 \right) (\bar{K}t/m)^{\beta m}
\nonumber\\
&\lesssim & \bar{K}^{-(1+\alpha)/2} ||f||_{\alpha} + 
\bar{K}^{1/2} \bar{K}_0  (\bar{K}t/m)^{\beta m}. \label{eq:critique} \EEA
The second term in (\ref{eq:critique}) is bounded by $c^{-1}\bar{K} (\bar{K}t/m)^{\beta m}$, in agreement with the
desired bound (\ref{eq:v-estimates}), but not the first one, which is bounded by $c^{-1}\bar{K} (\bar{K}t/(m-1))^{\beta(m-1)}$. 

\bigskip

In order to get an integrated bound of order $(\bar{K}t/m)^{\beta m}$ for the first term, we need a  refinement of Proposition \ref{prop:XJW}.  Fix $(t_1,x_1)\in Q^{(j)}$. We let (for $k\ge 0$ large enough so that $Q^{(j-k)}(t_1,x_1)\subset Q^{(j)} $)
\BEQ \tilde{v}^{(m)}(t,x):=v^{(m)}(t,x)+\int_{t}^{t_1} f(s,x_1)ds, \qquad (t,x)\in Q^{(j-k)}(t_1,x_1) \EEQ
so that $\tilde{v}^{(m)}$ satisfies the modified transport equation
\BEQ (\partial_{t'}-\Del+v^{(m-1)}\cdot\nabla)\tilde{v}^{(m)}(t,x)=\tilde{f}(t,x) \EEQ
with
\BEQ \tilde{f}(t,x):=f(t,x)-f(t,x_1).\EEQ
Note that $\nabla \tilde{v}^{(m)}=\nabla v^{(m)}, \nabla^2 \tilde{v}^{(m)}=\nabla v^{(m)}$. This introduces the following modifications. First, letting $\bar{B}_1^{(j-k)}:=\bar{B}(x_1,M^{(j-k)/2})$,
\BEQ ||\tilde{v}^{(m)}-v^{(m)}||_{\infty,Q^{(j-k)}(t_1,x_1)}\le \int_{t_1-M^j}^{t_1} ds ||f(s)||_{\infty,\bar{B}_1^{(j-k)}}
\le \bar{K}_0 (\bar{K}t/m)^{\beta m} \EEQ
as follows from (\ref{eq:3.8}), (\ref{eq:3.9}). Thus $||\tilde{v}^{(m)}||_{\infty,Q^{(j-k)}(t_1,x_1)}\lesssim 
\bar{K}_0 (\bar{K}t/m)^{\beta m}$ is bounded like $||v^{(m)}||_{\infty,Q^{(j-1)}}$.  Second (see (\ref{eq:ff})), 
$\tilde{f}(t,x)-\tilde{f}(t_1,x_1)=f(t,x)-f(t,x_1)$ involves values of $f$ {\em only} at time $t$. (Eventually this spares us having to bound inductively $\partial_t v^{(m)}$). 

\bigskip

We now go through
the proof of Proposition \ref{prop:XJW},  writing $\tilde{v}^{(m)}(t_1,x_1)$ as the sum of a series $\tilde{v}^{(m)}_{k_1+1}(t_1,x_1)+\sum_{k=k_1+1}^{\infty}
(\tilde{v}^{(m)}_{k+1}-\tilde{v}^{(m)}_k)(t_1,x_1)$, and bounding only $||\nabla \tilde{v}||_{\infty}=||\nabla v||_{\infty}$ and $||\nabla^2 \tilde{v}||_{\infty}=||\nabla^2 v||_{\infty}$.   Instead of (\ref{eq:estimatekk+1}), we get from the maximum principle
\BEQ \sup_{Q_1^{(j-1-k)}} |\tilde{v}^{(m)}_{k+1}-\tilde{v}^{(m)}_k|\lesssim M^{(j-k)(1+\alpha/2)} \left(
\fint_{t_1-M^{j-1-k}}^{t_1}  ds ||f(s)||_{\alpha,\bar{B}_1^{(j-1-k)}} + ||u^{(m-1)}||_{\alpha} \sup_{Q_1^{(j-1-k)}}\nabla \tilde{v}^{(m)} \right),\EEQ
where $\fint_{t-M^{j-1-k}}^t  (\ \cdot\ )ds:= M^{-(j-1-k)} \int_{t_1-M^{j-1-k}}^{t_1} (\ \cdot\ )ds$ is the average
over the time interval $[t_1-M^{j-1-k},t_1]$. We have proved above that $ ||f(s)||_{\alpha,\bar{B}_1^{(j)}}\lesssim
 c^{-1} \bar{K}^{(3+\alpha)/2} (\bar{K}s/(m-1))^{\beta(m-1)}$; thus (by explicit computation)
\BEA &&  \fint_{t_1-M^{j-1-k}}^{t_1}  ds ||f(s)||_{\alpha,\bar{B}_1^{(j-1-k)}}
\lesssim c^{-1} \bar{K}^{(3+\alpha)/2} \fint_{t-M^{j-1-k}}^t ds\, (\bar{K}s/(m-1))^{\beta(m-1)} \nonumber\\ && \qquad\qquad\qquad\qquad\qquad\qquad\qquad \equiv   c^{-1} \bar{K}^{(3+\alpha)/2}
(\bar{K}t/(m-1))^{\beta(m-1)} a_k,\EEA
with $a_k:=M^{k-j}t^{-\beta(m-1)} \frac{1}{\beta(m-1)+1} \left(t^{\beta(m-1)+1}-(t-M^{j-1-k})^{\beta(m-1)+1} \right)$.
Let $k_0:=\inf\{k\ge 0; M^{j-1-k}<t/m\}$; since $M^{j-1}\gtrsim t/m$ by hypothesis, $M^{j-1-k_0}\approx t/m$. For $k>k_0$, $a_k\approx 1$, as follows from Taylor's
formula; bounding all $a_k, k\ge 0$ by $1$
would yield the estimate (\ref{eq:critique}). However, for $k\le k_0$, $a_k\lesssim M^{k-j}\frac{t}{m}$, which is
a much better bound for $k_0-k$ large. Summarizing, the only change in the right-hand side of (\ref{eq:supnabla2}) is that $||f||_{\alpha}$ may be replaced by
\BEQ \sum_k  M^{-k\alpha/2} \fint_{t_1-M^{j-1-k}}^{t_1}  ds\,  ||f(s)||_{\alpha,\bar{B}_1^{(j-1-k)}}\lesssim
c^{-1}\bar{K}^{(3+\alpha)/2} (\bar{K}t/(m-1))^{\beta(m-1)} (A_1+A_2),\EEQ
where
\BEQ A_1:=\sum_{k\ge k_0} M^{-k\alpha/2}\lesssim (\bar{K}t/m)^{\alpha/2}\EEQ
and similarly
\BEQ A_2:=\sum_{k=0}^{k_0-1} M^{-k\alpha/2} M^{k-j}\frac{t}{m} \lesssim M^{k_0(1-\alpha/2)} (\bar{K}t/m)\approx
(\bar{K}t/m)^{\alpha/2}.\EEQ
All together, with respect to the rougher bound (\ref{eq:critique}), we have gained a small multiplicative
factor of order $A_1+A_2\lesssim (\bar{K}t/m)^{\beta}$, with $\beta:=\alpha/2$.  Thus
\BEA   ||\nabla v^{(m)}||_{\infty,Q^{(j-1)}}& \lesssim&  c^{-1}\bar{K} (\bar{K}t/(m-1))^{\beta(m-1)} \cdot
(\bar{K}t/m)^{\beta} +  c^{-1}\bar{K} (\bar{K}t/m)^{\beta m} \nonumber\\
&\lesssim & c^{-1}\bar{K} (\bar{K}t/m)^{\beta m}.\EEA

\end{itemize}
\hfill \eop


\section{H\"older estimates}


We prove in this section elementary H\"older estimates, together with a precise
form of the Schauder estimates which is crucial in the proof of Theorem \ref{th:main}
in section 3.

\begin{Definition}[H\"older semi-norms] 
Let $\gamma\in(0,1)$.
\begin{enumerate}
\item $f_0:\R^d\to\R$ is $\gamma$-H\"older continuous if $||f_0||_{\gamma}:=\sup_{x,x'\in\R^d} \frac{|f_0(x)-f_0(x')|}{|x-x'|^{\gamma}}<\infty$.
\item $f:\R_+\times\R^d\to\R$ is $\gamma$-H\"older continuous if $||f||_{\gamma}:=\sup_{(t,x),(t',x')\in\R_+\times\R^d} \frac{|f(t,x)-f(t',x')|}{|x-x'|^{\gamma}+|t-t'|^{\gamma/2}}<\infty$.
\end{enumerate}
\end{Definition}

In the denominator appearing in the definition of $||f||_{\gamma}$, we find a
power of the {\em parabolic distance}, $d_{par}((t,x),(t',x'))=|x-x'|+\sqrt{|t-t'|}$.
Note that $||\ ||_{\gamma}$ is only a semi-norm since $||1||_{\gamma}=0$. 
We also define H\"older semi-norms for functions restricted to $Q_0\subset\R_+\times\R^d$ or $Q\subset\R^d$ compact,
with the obvious definitions,
\BEQ ||f_0||_{\gamma,Q_0}:=\sup_{x,x'\in Q_0} \frac{|f_0(x)-f_0(x')|}{|x-x'|^{\gamma}}, \qquad
||f||_{\gamma,Q}:=\sup_{(t,x),(t',x')\in Q} \frac{|f(t,x)-f(t',x')|}{|x-x'|^{\gamma}+|t-t'|^{\gamma/2}}.\EEQ

{\em Remark.} For $f:\R_+\times\R^d\to\R$, we use in this article  either the
 parabolic H\"older semi-norm $||f||_{\alpha,Q}$ or the isotropic H\"older semi-norm $||f(t)||_{\alpha,Q_0}$ for
 $t\in\R_+$ fixed. The distinction is really important in the proof of Theorem \ref{th:v-estimates} (ii). Clearly,
 $||f(t)||_{\alpha,Q_0}\le ||f||_{\alpha,I\times Q_0}$ if $I$ is some time  interval containing $t$.

\begin{Lemma}[H\"older interpolation estimates] \label{lem:interpolation-Holder}
\begin{enumerate}
\item (on $\R^d$) Let $Q_0\subset\R$ be a convex set, and $u_0:Q_0\to\R$ such that
$||u_0||_{\infty,Q_0},||\nabla u_0||_{\infty,Q_0}<\infty$. Then
\BEQ ||u_0||_{\alpha,Q_0}\le ||u_0||^{1-\alpha}_{\infty,Q_0} ||\nabla u_0||^{\alpha}_{\infty,Q_0}, \qquad
\alpha\in(0,1).\EEQ
\item (on $\R_+\times\R^d$) 
Let $Q\subset\R_+\times\R^d$ be a convex set,  and $u:\R^d\to\R$ such that\\
$||u||_{\infty,Q},||\nabla u_0||_{\infty,Q},||\partial_t u_0||_{\infty,Q}<\infty$. Then
\BEQ ||u||_{\alpha,Q}\le 2\left( ||u||^{1-\alpha}_{\infty,Q} ||\nabla u||^{\alpha}_{\infty,Q} +
||u||_{\infty,Q}^{1-\alpha/2}||\partial_t u||_{\infty,Q}^{\alpha/2}\right), \qquad
\alpha\in(0,1).\EEQ
\end{enumerate}

\end{Lemma}

{\bf Proof.} (see \cite{Lie})  we prove (ii).  Let $X=(t,x)$ and $X'=(t',x')$ in $Q$, then
\BEA |u(X)-u(X')| &=& \left| \int_0^1 \frac{d}{d\tau} u((1-\tau)X+\tau X') d\tau \right| \nonumber\\
&\le & |t-t'|\ ||\partial_t u||_{\infty,Q} + |x-x'|\ ||\nabla u||_{\infty,Q} \le 2\max\left( |t-t'|\ ||\partial_t u||_{\infty,Q} , |x-x'|\ ||\nabla u||_{\infty,Q} \right). \nonumber\\ \EEA
On the other hand, $|u(X)-u(X')|\le 2||u||_{\infty}$. Hence
\BEQ |u(X)-u(X')|\le 2 \max\left( ||u||^{1-\alpha/2}_{\infty,Q} ||\partial_t u||^{\alpha/2}_{\infty,Q},
||u||^{1-\alpha}_{\infty,Q} ||\nabla u||^{\alpha}_{\infty,Q} \right).\EEQ
\hfill \eop

\begin{Lemma} \label{lem:Holder}
Let $u_0:\R^d\to\R$ be $\alpha$-H\"older. Then
\BEQ ||\nabla^{\kappa} (e^{t\Del}u_0)||_{\infty}\le C(d,\kappa,\alpha) t^{(\alpha-\kappa)/2} \,  ||u_0||_{\alpha}
\qquad (\kappa\ge 1);
\label{eq:bound-nabla2} \EEQ
\BEQ  ||\nabla^2 (e^{t\Del}u_0)||_{\gamma}\le C'(d,\gamma,\alpha) t^{-1+(\alpha-\gamma)/2} \, ||u_0||_{\alpha}
\qquad (\gamma\in (0,1)); \label{eq:bound-nabla2alpha}\EEQ
\BEQ ||e^{t\Del}u_0-e^{t'\Del}u_0||_{\infty}\le C''(d,\alpha) (t-t')^{\alpha/2}\, 
||u_0||_{\alpha} \qquad (\alpha\in(0,1),\ t>t'>0).\label{eq:bound-nabla2t} \EEQ
\end{Lemma}

{\bf Proof.} 
 (\ref{eq:bound-nabla2alpha}) follows by Lemma \ref{lem:interpolation-Holder} from the bounds (\ref{eq:bound-nabla2})
 with $\kappa=2,3$. Thus let us first prove (\ref{eq:bound-nabla2}). The regularizing operator $e^{t\Del}$ is
 defined by convolution with respect to the heat kernel $p_t$. By translation invariance, it is enough to bound the
 quantity $I(\eps):=\nabla^{\kappa-1}(e^{t\Del}u_0)(0)-
\nabla^{\kappa-1}(e^{t\Del}u_0)(\eps)$ in the limit $\eps\to 0$. The quantities in (\ref{eq:bound-nabla2})
are invariant through the substitution $u_0\to u_0-u_0(0)$, so we assume that $u_0(0)=0$.  We may also assume $|\eps|\ll \sqrt{t}$. Let
$A:=\eps^{\beta}t^{(1-\beta)/2}$ with $\beta=(1-\alpha)/d$; note that $|\eps|\ll A\ll \sqrt{t}$.  We split the integral into three parts, $I(\eps)=I_1(\eps)+I_2(\eps)+I_3(\eps)$,  with
\BEQ I_1(\eps):=\int_{|x|< A} dx\, \nabla^{\kappa-1} p_t(x) (u_0(x)-u_0(x+\eps)), \ \  I_2(\eps):=\int_{|x|> A} dx \, (\nabla^{\kappa-1} p_t(x)-\nabla^{\kappa-1} p_t(x+\eps)) (u_0(x)-u_0(0))\EEQ
\BEQ I_3(\eps)=\left( \int_{|x|>A} dx -\int_{|x-\eps|>A} dx \right) \nabla^{\kappa-1} p_t(x+\eps) (u_0(x)-u_0(0)).\EEQ
We  use
$|u_0(x)-u_0(x+\eps)|\le ||u_0||_{\alpha}\, |\eps|^{\alpha}$ in the first integral, and get
\BEQ I_1(\eps)\lesssim  ||u_0||_{\alpha}  A^d t^{-(\kappa+d-1)/2} |\eps|^{\alpha}= ||u_0||_{\alpha} t^{(\alpha-\kappa)/2}|\eps|. \label{eq:I1eps} \EEQ
For the second integral, we use 
$|\nabla^{\kappa-1} p_t(x)-\nabla^{\kappa-1} p_t(x+\eps)|\lesssim \frac{|\eps|}{t^{\kappa/2}} p_t(x)$ and 
$|u_0(x)-u_0(0)|\le ||u_0||_{\alpha} |x|^{\alpha}$, yielding the same estimate. Finally, the integration
volume in the third integral is $O(A^{d-1}|\eps|)$, hence $I_3(\eps)\lesssim  ||u_0||_{\alpha} A^{d-1}|\eps| t^{-(\kappa-1)/2} A^{\alpha}
\lesssim ||u_0||_{\alpha} A^d t^{-(\kappa+d-1)/2} |\eps|^{\alpha} \ \cdot\ (|\eps|/A)^{1-\alpha}$ is negligible with respect to
the first integral (compare with (\ref{eq:I1eps})).
 Taking $\eps\to 0$, this
gives the desired bound for $ ||\nabla^{\kappa} (e^{t\Del}u_0)||_{\infty}$. 

\medskip
Finally, (\ref{eq:bound-nabla2t}) may be obtained through the use of the fractional
derivative $$|\nabla|^{\alpha}:u_0\mapsto \left(|\nabla|^{\alpha}u_0:x\mapsto \int d\xi dy
|\xi|^{\alpha} e^{\II(x-y)\xi} u_0(y) \right),$$ namely,
\BEA |(e^{t\Del}u_0-e^{t'\Del}u_0)(x)|&=& \left| \int_{t'}^t ds \int dy\,  \partial_s p_s(x-y) u_0(y) \right| = \left| \int_{t'}^t ds \int dy\,  \Del p_s(x-y) u_0(y) \right| \nonumber\\
&\lesssim & \int_{t'}^t ds \int dy\, |\,  |\nabla|^{2-\alpha/2} p_s(x-y)| \ \ | \, |\nabla|^{\alpha}u_0(y)|
\lesssim (t^{\alpha/2}-(t')^{\alpha/2}) \, ||u_0||_{\alpha} \nonumber\\
& \lesssim & (t-t')^{\alpha/2}
||u_0||_{\alpha}. \EEA 
 \hfill \eop

\begin{Corollary} \label{cor:Holder}
Let $g:[0,t]\times\R^d\to\R$ be a continuous function such that $(g_s)_{s\in[0,t]}$ are
uniformly $\alpha$-H\"older, and $\gamma<\alpha$. Then $s\mapsto ||\nabla^2 (e^{(t-s)\Del}g_s)||_{\gamma}$ is $L^1_{loc}$ and, for $0<t'<t$,
\BEQ \int_{t'}^t ds \, ||\nabla^2 (e^{(t-s)\Del}g_s)||_{\gamma}  \le C''(d,\gamma,\alpha) (t-t')^{(\alpha-\gamma)/2} \sup_{s\in[t',t]} ||g_s||_{\alpha} .\EEQ
\end{Corollary}

\bigskip

We now turn to  our Schauder estimates. The multi-scale proof of the Proposition below is  inspired by Wang \cite{XJW}. We fix  a constant $M>1$, e.g. $M=2$ for
a dyadic scale decomposition. 

\begin{Definition}[parabolic balls] 
Let $(t_0,x_0)\in\R\times\R^d$ and $j\in\Z$. Then the scale $j$ parabolic ball issued 
from $(t_0,x_0)$ is the closed subset $Q^{(j)}(t_0,x_0):=\{(t,x)\in\R\times\R^d;\ t_0-M^{j}\le t\le t_0, x\in \bar{B}(x_0,M^{j/2})\}$.
\end{Definition}

The set $\{(t,x)\ |\ t\le t_0, d_{par}((t,x),(t_0,x_0))\le M^{j/2}\}$ is comparable
to $Q^{(j)}(t_0,x_0)$, in the sense that there exist $\del k_0,\del k_1\ge 0$ such that
$Q^{(j)}(t_0,x_0)\subset \{(t,x)\ |\ t\le t_0, d_{par}((t,x),(t_0,x_0))\le M^{(j+\del k_0)/2}\}
\subset Q^{(j+\del k_0+\del k_1)/2}(t_0,x_0)$ (one may actually choose $\del k_1=0$), which is why $Q^{(j)}(t_0,x_0)$ is called a 'ball';
but mind the causality condition $t\le t_0$. In the sequel we let $\del k=\del k(M)$ be some large
enough integer, depending only on $M$, used in several occasions to make different parabolic balls fit exactly into each other. The main property of parabolic balls
in our context is the simple scaling property for locally bounded solutions $u$ of the heat equation
$(\partial_t-\Del)u=0$: for all $\kappa=(\kappa_1,\ldots,\kappa_d)$, $\kappa_1,\ldots,\kappa_d\ge 0$,  $|\nabla^{\kappa} u(t_0,x_0)|\lesssim (M^{-j/2})^{|\kappa|} \sup_{\partial_{par}Q^{(j)}(t_0,x_0)} |u|$ ($|\kappa|=\kappa_1+\ldots+\kappa_d$), where $\partial_{par}Q^{(j)}(t_0,x_0):= \left(\{t_0-M^j\}\times \bar{B}(x_0,M^{j/2})\right)\cup
\left( [t_0-M^j,t_0)\times \partial B(x_0,M^{j/2}) \right)$ is the {\em parabolic boundary} of $Q^{(j)}(t_0,x_0)$.  From this we simply deduce the following: 
let 
\BEQ Q^{(j)}_{(k)}(t_0,x_0):=\{(t,x)\in Q^{(j)}(t_0,x_0)\ |\ d_{par}((t,x),\partial_{par}Q^{(j)}(t_0,x_0))\ge M^k\} \qquad  (k\le j), \EEQ
 then $\sup_{Q^{(j)}_{(k)}(t_0,x_0)} 
|\nabla^{\kappa} u|\lesssim (M^{-k/2})^{|\kappa|} \sup_{Q^{(j)}(t_0,x_0)} |u|$, which is a quantitative version of the well-known regularizing
property of the heat equation: if $u$ is bounded on some $j$ scale parabolic ball $Q^{(j)}$, then $\nabla^{\kappa}u$ is bounded {\em away} from the parabolic boundary
of $Q^{(j)}$.  In particular, since $Q^{(j-1)}(t_0,x_0)\subset Q^{(j)}_{(j-\del k)}(t_0,x_0)$, one has: $\sup_{Q^{(j-1)}(t_0,x_0)} 
|\nabla^{\kappa} u|\lesssim (M^{-j/2})^{|\kappa|} \sup_{Q^{(j)}(t_0,x_0)} |u|$.

\medskip

\begin{Proposition}[Schauder estimates]  \label{prop:XJW}
Let $v$ solve the linear parabolic PDE 
\BEQ (\partial_t-\Del+a(t,x))u(t,x)=b(t,x)\cdot
\nabla u(t,x) + f(t,x) \EEQ
 on the parabolic ball 
$Q^{(j)}:=Q^{(j)}(t_0,x_0)$. Assume: $u$ is bounded; $a\ge 0$;
 
\BEQ ||f||_{\alpha}:=||f||_{\alpha,Q^{(j)}}:=\sup_{(t,x),(t',x')\in Q^{(j)}} \frac{|f(t,x)-f(t',x')|}{|x-x'|^{\alpha}+|t-t'|^{\alpha/2}}<\infty \EEQ for some $\alpha\in(0,1)$,
and similarly $||a||_{\alpha},||b||_{\alpha}<\infty$.
 Then 
\BEQ 
\sup_{Q^{(j-1)}} |\nabla u| \lesssim M^{j/2} R_b^{-1} \left\{ 
 M^{j\alpha/2} ||f||_{\alpha}+ \left(M^{j\alpha} R_b^{-1} ||b||_{\alpha}^2+
 M^{j\alpha/2} ||a||_{\alpha}+  M^{-j} \right) \sup_{Q^{(j)}} |u| \right\}, \label{eq:3.46} \EEQ
 
\BEA && 
 ||\nabla u||_{\alpha,Q^{(j-1)}} \lesssim  M^{-j\alpha/2} R_b^{-(1+\alpha)/2} \  \left\{ 
 M^{j(1+\alpha)/2}  ||f||_{\alpha}  \right. \nonumber\\
 && \qquad \qquad   \left.  + \left(M^{j(1+\alpha+\alpha^2)/2\alpha}  R_b^{-\half(1+\alpha)/\alpha} ||b||_{\alpha}^{(1+\alpha)/\alpha} + M^{j(1+\alpha)/2} ||a||_{\alpha}+ M^{-j/2}\right) \sup_{Q^{(j)}} |u| \right\}, 
 \label{eq:3.47} \EEA
 
\BEQ 
\sup_{Q^{(j-1)}} |\partial_t u|, \sup_{Q^{(j-1)}} |\nabla^2 u| \lesssim  R_b^{-1} \left\{ 
 M^{j\alpha/2}  ||f||_{\alpha}+ \left(M^{j\alpha} R_b^{-1} ||b||_{\alpha}^2+
 M^{j\alpha/2} ||a||_{\alpha}+ M^{-j} \right) \sup_{Q^{(j)}} |u| \right\}, 
\label{eq:3.48} \EEQ

and for every $\alpha'>\alpha$,
\BEA && 
||\partial_t u||_{\alpha,Q^{(j-1)}}, ||\nabla^2 u||_{\alpha,Q^{(j-1)}} \lesssim M^{-j\alpha/2} R_b^{-(1+\alpha'/2)}\  \left\{ 
 M^{j\alpha/2}  ||f||_{\alpha}  \right. \nonumber\\
 && \qquad \qquad \qquad  \left.  + \left(M^{j\alpha/2}  R_b^{-\half(2+\alpha')/(1+\alpha)} ||b||_{\alpha}^{(2+\alpha)/(1+\alpha)} + M^{j\alpha/2} ||a||_{\alpha}+ M^{-j}\right) \sup_{Q^{(j)}} |u| \right\}, 
 \label{eq:3.49}
 \EEA
 where $R_b:=\left(1+M^{j/2}|b(t_0,x_0)|\right)^{-1}$.
\end{Proposition} 

\medskip

{\em Remark:}
Removing the condition $a\ge 0$, we would get the same estimates, multiplied by
$e^{M^j \sup_{Q^{(j)}} (-a)}$. 

\medskip

{\bf Proof.} 
Let $\tilde{u}(\tilde{t},\tilde{x}):=u(M^j \tilde{t},M^{j/2}\tilde{x})$, $\tilde{b}(\tilde{t},\tilde{x}):=M^{j/2}
b(M^j \tilde{t},M^{j/2}\tilde{x})$, $\tilde{f}(\tilde{t},\tilde{x}):=M^j f(M^j\tilde{t},M^{j/2}\tilde{x})$, $\tilde{a}(\tilde{t},\tilde{x}):=M^j a(M^j \tilde{t}, 
M^{j/2} \tilde{x}).$ 
Then the PDE $(\partial_t-\Del+a)u=b\cdot\nabla u+f$ on $Q^{(j)}$ reduces to an equivalent PDE,
$(\partial_{\tilde{t}}-\tilde{\Del}+\tilde{a})\tilde{u}=\tilde{b}\cdot\tilde{\nabla}\tilde{u}+\tilde{f}$ on a parabolic ball
$\tilde{Q}$ of size unity. {\em Assume } (leaving out for sake of conciseness the powers of $R_b=
(1+|\tilde{b}(t_0,x_0)|)^{-1}$) that we have proved an inequality of the type
\BEQ \sup_{\tilde{Q}^{(-1)}} |\tilde{\nabla}^{\kappa} \tilde{u}|\lesssim  \left(
||\tilde{f}||_{\alpha}+(||\tilde{b}||_{\alpha}^{\beta}+||\tilde{a}||_{\alpha}+1)\sup_{\tilde{Q}}|\tilde{u}| \right), \qquad \qquad {\mathrm{resp.}} \EEQ
\BEQ ||\tilde{\nabla}^{\kappa} \tilde{u}||_{\alpha,\tilde{Q}^{(-1)}}\lesssim  \left(
||\tilde{f}||_{\alpha}+(||\tilde{b}||_{\alpha}^{\beta}+||\tilde{a}||_{\alpha}+1)\sup_{\tilde{Q}}|\tilde{u}| \right). \EEQ
By rescaling, we get
\BEQ \sup_{Q^{(j-1)}} |\nabla^{\kappa} u|\lesssim (M^{-j/2})^{\kappa}   \left( M^{j(1+\alpha/2)}
||f||_{\alpha}+\left( ((M^{j/2})^{1+\alpha} ||b||_{\alpha})^{\beta}+ (M^j)^{1+\alpha/2} 
||a||_{\alpha}+1 \right)\sup_{Q^{(j)}}|u| \right),  \EEQ
\BEQ ||\tilde{\nabla}^{\kappa} u||_{\alpha,Q^{(j-1)}}\lesssim (M^{-j/2})^{\kappa+\alpha}  \left(
M^{j(1+\alpha/2)} ||f||_{\alpha}+\left( ((M^{j/2})^{1+\alpha} ||b||_{\alpha})^{\beta}+
(M^j)^{1+\alpha/2} ||a||_{\alpha} +1 \right)
\sup_{Q^{(j)}}|u| \right). \EEQ

This gives the correct scaling factors in (\ref{eq:3.46},\ref{eq:3.47},\ref{eq:3.48},\ref{eq:3.49}). Thus we may 
assume that $j=0$. In the sequel we write for short $||\ \cdot \ ||_{\alpha}$ instead of $||\ \cdot\ ||_{\alpha,Q^{(0)}}$
and $||\ \cdot \ ||_{\infty}$ instead of $\sup_{Q^{(0)}} |\ \cdot\ |$. 

\medskip

The general principle underlying the proof
of the Schauder estimates in 
 \cite{XJW} is the following. Let $(t_1,x_1)\in Q^{(0)}_{(-k_1)}$. One rewrites $u(t_1,x_1)$ as the sum of the series $u(t_1,x_1)=u_{k_1+1}(t_1,x_1)+\sum_{k=k_1+1}^{+\infty} (u_{k+1}(t_1,x_1)-u_k(t_1,x_1))$,
 where $u_k$, $k\ge k_1+1$ is the solution on $Q_1^{(-k)}:=Q^{(-k)}(t_1,x_1)$ of the 'frozen' PDE 
 \BEQ (\partial_t-\Del+a(t_1,x_1))u_k(t,x)=b(t_1,x_1)\cdot 
 \nabla u_k(t,x)+f(t_1,x_1) \label{eq:frozen} \EEQ
 with initial-boundary condition $u_k\big|_{\partial_{par}Q_1^{(-k)}}=u\big|_{\partial_{par}Q_1^{(-k)}}$. We split the
 proof into several steps.
 
\begin{itemize}
\item[(i)] (estimates for $|u_{k+1}-u_k|$) One first remarks that $u_k-u$, $k\ge k_1+1$ solves on $Q_1^{(-k)}$ the heat equation
\BEQ (\partial_t-\Del+a(t_1,x_1)-b(t_1,x_1)\cdot\nabla)(u_k-u)=(b(t_1,x_1)-b)\cdot \nabla u+(f(t_1,x_1)-f)-(a(t_1,x_1)-a)u \label{eq:ff} \EEQ
with zero initial-boundary condition  $(u_k-u)\big|_{\partial_{par}
Q_1^{(-k)}}=0$,
implying by the maximum principle
\BEQ \sup_{Q_1^{(-k-1)}} |u_{k+1}-u_k|\le \sup_{Q_1^{(-k-1)}} |u_{k+1}-u| +  \sup_{Q_1^{(-k)}} |u_{k}-u|\lesssim
 M^{-k(1+\alpha/2)} \left( ||f||_{\alpha} + ||a||_{\alpha} ||u||_{\infty} + ||b||_{\alpha} \sup_{Q_1^{(-k)}} |\nabla u| \right).
 \label{eq:estimatekk+1} \EEQ
 
\item[(ii)] (estimates for higher-order derivatives of $u_{k_1+1}$)
Recall $u_{k_1+1}$ is a solution of the heat equation $(\partial_t-\Del-b(t_1,x_1)\cdot\nabla)u_{k_1+1}=f(t_1,x_1)$ with initial-boundary
condition $u_{k_1+1}\big|_{\partial_{par}Q_1^{-(k_1+1)}}=u\big|_{\partial_{par}Q_1^{-(k_1+1)}}$.

Assume first $|b(t_1,x_1)|\lesssim 1$. As follows from standard estimates recalled before
the proposition,
\BEQ ||\nabla u_{k_1+1}||_{\alpha,Q_1^{-(k_1+2)}} \lesssim (M^{k_1/2})^{1+\alpha} ||u||_{\infty} , \ \  \sup_{Q_1^{-(k_1+2)}} |\partial_t u_{k_1+1}|, \sup_{Q_1^{-(k_1+2)}} |\nabla^2 u_{k_1+1}| \lesssim
 M^{k_1} ||u||_{\infty} , \EEQ
 \BEQ \qquad\qquad\qquad\qquad   ||\nabla^2 u_{k_1+1}||_{\alpha,Q_1^{-(k_1+2)}} \lesssim (M^{k_1})^{1+\alpha/2} ||u||_{\infty}.\EEQ

 If $|b(t_0,x_0)|\gg 1$, then one makes the Galilean transformation
$x\mapsto x-b(t_0,x_0)t$ to get rid of the drift, after which the boundary of $Q_1^{-(k_1+1)}$ lies at distance $R=O(M^{-k_1/2}/|b(t_0,x_0)|)$ instead of $O(M^{-k_1/2})$ of $(t_1,x_1)$; thus, in general,

\BEQ ||\nabla u_{k_1+1}||_{\alpha,Q_1^{-(k_1+2)}} \lesssim  R_b^{-(1+\alpha)/2} (M^{k_1/2})^{1+\alpha} ||u||_{\infty} , \ \  \sup_{Q_1^{-(k_1+2)}} |\partial_t u_{k_1+1}|, \sup_{Q_1^{-(k_1+2)}} |\nabla^2 u_{k_1+1}| \lesssim R_b^{-1}
 M^{k_1} ||u||_{\infty} , \EEQ
 \BEQ \qquad\qquad\qquad\qquad   ||\nabla^2 u_{k_1+1}||_{\alpha,Q_1^{-(k_1+2)}} \lesssim
 R_b^{-(1+\alpha/2)} (M^{k_1})^{1+\alpha/2} ||u||_{\infty}.\EEQ

\item[(iii)] (estimates for higher-order derivatives of $u_{k+1}-u_k$)
Similarly to (ii), we note that $u_{k+1}-u_k$ is a solution on $Q_1^{(-k-1)}$ of the heat equation
$(\partial_t-\Del+a(t_1,x_1)-b(t_1,x_1)\cdot\nabla)(u_{k+1}-u_k)=0$.
Thus
\BEA && \sup_{Q_1^{(-k-2)}} |\partial_t (u_{k+1}-u_k)|, \sup_{Q_1^{(-k-2)}} |\nabla^2 (u_{k+1}-u_k)| \lesssim M^k R_b^{-1} \sup_{Q_1^{(-k-1)}} |u_{k+1}-u_k|,
\label{eq:Schauder-iii}  \\
&& \qquad \qquad 
 ||\nabla^2 (u_{k+1}-u_k)||_{\alpha',Q_1^{(-k-2)}} \lesssim (M^k)^{1+\alpha'/2} R_b^{-(1+\alpha'/2)} 
\sup_{Q_1^{(-k-1)}} |u_{k+1}-u_k| \label{eq:Schauder-iii-alpha} \EEA
is bounded using (i) in terms of $R_b$, $||b||_{\alpha},||f||_{\alpha}$ and $\sup_{Q_1^{(-k)}} |\nabla u|.$

\item[(iv)] (Schauder estimates for higher-order derivatives of $u$)
Summing up the estimates in (i), (ii), (iii), and noting that $\cdots \subset 
Q_1^{(-k_1-2)}\subset Q_1^{(-k_1-1)}\subset Q^{(0)}_{(-k_1-\del k)}$ for $\del k=\del k(M)$ large enough, one obtains
\BEQ M^{-k_1} \sup_{Q_{(-k_1)}^{(0)}} |\partial_t u|, M^{-k_1} \sup_{Q^{(0)}_{(-k_1)}} |\nabla^2 u|\lesssim R_b^{-1} \left\{ (M^{-k_1})^{1+\alpha/2} \left( ||f||_{\alpha}+
||a||_{\alpha} ||u||_{\infty}+ ||b||_{\alpha}\sup_{Q^{(0)}_{(-k_1-\del k})} |\nabla u|\right)+
||u||_{\infty} \right\}. \label{eq:supnabla2} \EEQ

By interpolation (see immediately thereafter),  $\sup_{Q^{(0)}_{(-k_1-\del k)}} |\nabla u|$ is bounded
in terms of  $||u||_{\infty}$ and $\sup_{Q^{(0)}_{(-k_1-\del k)}} |\nabla^2 u|$. Thus
in principle (\ref{eq:supnabla2}) gives a bound for $\nabla^2 u$. {\em However}, since
$Q^{(0)}_{(-k_1-\del k)}\supsetneq Q^{(0)}_{(-k_1)}$, one  {\em cannot} fix $k_1$. Instead we shall bound $\sup_{k_1} M^{-k_1} \sup_{Q^{(0)}_{(-k_1)}} |\nabla^2 u|$, and similarly
for the different gradient/H\"older norms considered in the Proposition. This
explains {\em why} ultimately we must consider the values of $\nabla u$, $\nabla^2 u$
on the whole parabolic ball $Q^{(0)}$, not only on the subset $Q^{(-1)}$ where our
results are stated. 

Now 
\BEQ \sup_{Q^{(0)}_{(-k_1-\del k)}} |\nabla u|\lesssim
 \left(\sup_{Q^{(0)}_{(-k_1-\del k)} }
|\nabla^2 u| \right)^{1/2}  \left( ||u||_{\infty} \right)^{1/2} \lesssim \eps^2 \sup_{Q^{(0)}_{(-k_1-\del k)} }
|\nabla^2 u| +  \eps^{-2} ||u||_{\infty}  \label{eq:interpol-nabla} \EEQ
 for every $\eps>0$.  Hence (using (\ref{eq:supnabla2})), choosing $\eps^2\approx
 R_b/||b||_{\alpha}$, one gets

\BEQ \sup_{k_1\ge 0}  M^{-k_1} \sup_{Q_{(-k_1)}^{(0)}} |\nabla^2 u|   \lesssim
R_b^{-1} \left\{ (M^{-k_1})^{1+\alpha/2} \left( ||f||_{\alpha}+ (||a||_{\alpha}+R_b^{-1} ||b||_{\alpha}^2) ||u||_{\infty} \right)+ ||u||_{\infty}  \right\} ,\EEQ
implying  in particular the bound (\ref{eq:3.48}) for $\nabla^2 u$, from which 
(\ref{eq:interpol-nabla},\ref{eq:supnabla2}) yields the bound (\ref{eq:3.48}) for $\partial_t u$.

Using the estimates (\ref{eq:3.48}) and (\ref{eq:interpol-nabla}) with $\eps=1$
yields also the gradient bound (\ref{eq:3.46}).

\item[(v)] (Schauder estimates for H\"older norms)

Let us now bound $||\nabla^2 u||_{\alpha,Q^{(0)}_{-(k_1-1)}}\approx \sup_{(t_1,x_1),(t_2,x_2)\in Q^{(0)}_{-(k_1-1)}} \frac{|\nabla^2 u(t_2,x_1)-\nabla^2 u(t_2,x_2)|}{
d_{par}((t_1,x_1),(t_2,x_2))^{\alpha}}$ or equivalently
$||\partial_t u||_{\alpha,Q^{(0)}_{(-k_1-1)}}$.  Assume e.g. $t_1\ge t_2$, and
$(t_2,x_2)\in Q^{(-k_2)}(t_1,x_1)$, $k_2\ge k_1+1$, with $d_{par}((t_1,x_1),(t_2,x_2))
\approx M^{-k_2/2}$. The hypothesis $k_2\ge k_1+1$ excludes the case where 
$d_{par}((t_1,x_1),(t_2,x_2))$ is comparable to $M^{-k_1/2}$, a case which is not needed
since it is already covered by the estimates proved in (iv). 
Then  $|\nabla^2 u(t,x)-\nabla^2 u(t',x')|\le I_1+I_2+I_3+I_4$, with (using
(\ref{eq:Schauder-iii-alpha}) for $I_1,I_2$ and (\ref{eq:Schauder-iii}) for $I_3,I_4$)
\BEQ I_1=|\nabla^2 u_{k_1}(t_1,x_1)-\nabla^2 u_{k_1}(t_2,x_2)|\lesssim (M^{k_1})^{1+\alpha/2}
R_b^{-(1+\alpha/2)}  ||u||_{\infty} d_{par}(t_1,x_1;t_2,x_2)^{\alpha};\EEQ
\BEA  I_2 &=& \sum_{k=k_1}^{k_2-1} |\nabla^2 (u_{k+1}-u_{k})(t_1,x_1)-\nabla^2 (u_{k+1}-u_{k})(t_2,x_2)|
\nonumber\\  &\lesssim &  R_b^{-(1+\alpha'/2)}
d_{par}(t_1,x_1;t_2,x_2)^{\alpha'} \left(\sum_{k=k_1}^{k_2-1} (M^{k/2})^{\alpha'-\alpha}\right)  \left(  ||f||_{\alpha} + ||a||_{\alpha} ||u||_{\infty} +||b||_{\alpha}\sup_{Q_1^{(-k)}}|\nabla u|\right)  \nonumber\\
&\lesssim& d_{par}(t_1,x_1;t_2,x_2)^{\alpha} R_b^{-(1+\alpha'/2)}   \left(  ||f||_{\alpha}
+ ||a||_{\alpha} ||u||_{\infty} +||b||_{\alpha}\sup_{Q_1^{(-k_1)}}|\nabla u|\right); \label{eq:4.19} \EEA
and \BEQ I_3:=
\sum_{k\ge k_2} |\nabla^2 (u_{k+1}-u_k)(t_1,x_1)|, \  I_4:=\sum_{k\ge k_0} |\nabla^2 (u_{k+1}-u_k)(t_2,x_2)| \EEQ 
are
\BEQ \lesssim  d_{par}(t_1,x_1;t_2,x_2)^{\alpha} R_b^{-1}\left(||f||_{\alpha} + ||a||_{\alpha} ||u||_{\infty} + ||b||_{\alpha} \sup_{Q_1^{(-k_2)}} |\nabla u| \right)\EEQ

Hence 
\BEA && (M^{-k_1})^{1+\alpha/2} ||\partial_t u||_{\alpha,Q^{(0)}_{-(k_1-1)}},  (M^{-k_1})^{1+\alpha/2} ||\nabla^2 u||_{\alpha,Q^{(0)}_{-(k_1-1)}} \lesssim R_b^{-(1+\alpha'/2)} \cdot \nonumber\\
&& \qquad \cdot\  \ 
 \left\{ (M^{-k_1})^{1+\alpha/2} \left(  ||f||_{\alpha}+||b||_{\alpha}\sup_{Q^{(0)}_{(-k_1-\del k)}} |\nabla u|+ ||a||_{\alpha} ||u||_{\infty} \right) + 
||u||_{\infty} \right),  \label{eq:4.21} \EEA
compare with (\ref{eq:supnabla2}).

By standard H\"older interpolation inequalities \cite{Lie},
\BEQ \sup_{Q^{(0)}_{(-k_1-\del k)}}|\nabla u|\lesssim  ||\nabla^2 u||_{\alpha,Q^{(0)}_{(-k_1-\del k)}}^{1/(2+\alpha)} \ (\sup_{Q^{(0)}_{(-k_1-\del k)}} |u|)^{(1+\alpha)/(2+\alpha)} \lesssim 
\eps^{2+\alpha} ||\nabla^2 u||_{\alpha,Q^{(0)}_{(-k_1-\del k)}}+\eps^{-(2+\alpha)/(1+\alpha)} ||u||_{\infty} 
\label{eq:interpol-alpha} \EEQ
for every $\eps>0$. Choosing $\eps^{2+\alpha}\approx R_b^{1+\alpha'/2}/||b||_{\alpha}$ yields as in (iv) a bound for\\ $\sup_{k_1\ge 0} (M^{-k_1})^{1+\alpha/2} ||\nabla^2 u||_{\alpha,Q^{(0)}_{-(k_1-1)}}$, from which one deduces in particular
(\ref{eq:3.49}).

\medskip

In order to obtain the bound (\ref{eq:3.47}) for $||\nabla u||_{\alpha,Q^{(-1)}}$, we proceed initially in the same way, with the only difference that one may take $\alpha'=\alpha$
in (\ref{eq:4.19}) since one gets a series $\sum_{k=k_1}^{k_2-1} M^{-k/2}$  of order $O(1)$. Thus (\ref{eq:4.21}) becomes
\BEQ  (M^{-k_1/2})^{1+\alpha} ||\nabla u||_{\alpha,Q^{(0)}_{-(k_1-1)}} \lesssim R_b^{-(1+\alpha)/2}
 \left\{ (M^{-k_1/2})^{1+\alpha} \left( ||f||_{\alpha}+||b||_{\alpha}\sup_{Q^{(0)}_{(-k_1-\del k)}} |\nabla u|+ ||a||_{\alpha} ||u||_{\infty} \right) + 
||u||_{\infty} \right\}. \label{eq:4.21bis} \EEQ
One now uses H\"older interpolation inequalities to bound $\nabla u$ in terms of $||u||_{\infty}$
and $\nabla^2 u$. Instead of (\ref{eq:interpol-alpha}), one has here
\BEQ \sup_{Q^{(0)}_{(-k_1-\del k)}}|\nabla u|\lesssim  ||\nabla u||_{\alpha,Q^{(0)}_{(-k_1-\del k)}}^{1/(1+\alpha)} \ (\sup_{Q^{(0)}_{(-k_1-\del k)}} |u|)^{\alpha/(1+\alpha)} \lesssim 
\eps^{1+\alpha} ||\nabla u||_{\alpha,Q^{(0)}_{(-k_1-\del k)}}+\eps^{-(1+\alpha)/\alpha} ||u||_{\infty} 
\label{eq:interpol-alpha} \EEQ
for every $\eps>0$. Choosing $\eps^{1+\alpha}\approx R_b^{(1+\alpha)/2}/||b||_{\alpha}$ yields as in (iv) a bound for\\ $\sup_{k_1\ge 0} (M^{-k_1})^{(1+\alpha)/2} ||\nabla u||_{\alpha,Q^{(0)}_{-(k_1-1)}}$, from which one deduces in particular
(\ref{eq:3.47}).

\end{itemize}

\hfill\eop



\end{document}